\documentclass{amsart}
\usepackage[all]{xy}
\usepackage{latexsym}
\usepackage{amssymb}
\usepackage{amsmath}
\usepackage{amsthm}
\usepackage{mathabx}
\usepackage{amscd}
\usepackage{fancyhdr}
\usepackage{fullpage}
\usepackage{mathrsfs}
\xyoption{arc}

 \usepackage{tikz}
\usetikzlibrary{matrix}
\usetikzlibrary{shapes}
\usetikzlibrary{arrows}
\usetikzlibrary{calc,3d}
\usetikzlibrary{decorations,decorations.pathmorphing,decorations.pathreplacing,decorations.markings}
\usetikzlibrary{through}

\tikzset{snake arrow/.style=
{->,
decorate,
decoration={snake,amplitude=.4mm,segment length=2mm,post length=1mm}},
}

\tikzset{ext/.style={circle, draw,inner sep=1pt},int/.style={circle,draw,fill,inner sep=1.4pt},nil/.style={inner sep=1pt}}
\tikzset{cy/.style={circle,draw,fill,inner sep=2pt},scy/.style={circle,draw,inner sep=2pt},scyx/.style={draw,cross out,inner sep=2pt},scyt/.style={draw,regular polygon,regular polygon sides=3,inner sep=0.95pt}}
\tikzset{exte/.style={circle, draw,inner sep=3pt},inte/.style={circle,draw,fill,inner sep=3pt}}
\tikzset{diagram/.style={matrix of math nodes, row sep=3em, column sep=2.5em, text height=1.5ex, text depth=0.25ex}}
\tikzset{diagram2/.style={matrix of math nodes, row sep=0.5em, column sep=0.5em, text height=1.5ex, text depth=0.25ex}}

\tikzset{
  rightblue/.style={
    decoration={markings,mark=at position .8 with {\arrow[scale=1.2,blue]{latex}}},
    postaction={decorate},
    shorten >=0.4pt}}
\tikzset{
  leftblue/.style={
    decoration={markings,mark=at position .55 with {\arrowreversed[scale=1.2,blue]{latex}}},
    postaction={decorate},
    shorten >=0.4pt}}
\tikzset{
  rightred/.style={
    decoration={markings,mark=at position .45 with {\arrow[scale=1.2,red]{latex}}},
    postaction={decorate},
    shorten >=0.4pt}}
\tikzset{
  leftred/.style={
    decoration={markings,mark=at position .2 with {\arrowreversed[scale=1.2,red]{latex}}},
    postaction={decorate},
    shorten >=0.4pt}}

  \def\id{{\mbox{1 \hskip -7pt 1}}}
\newcommand{\sgn}{{\mathit s  \mathit g\mathit  n}}
 \newcommand{\lon}{\longrightarrow}
 \newcommand{\bu}{\bullet}
 
 \newcommand{\rar}{\rightarrow}
 \newcommand{\hook}{\hookrightarrow}

\newcommand{\p}{{\partial}}

\newcommand{\Der}{\mathrm{Der}}

\newcommand{\wh}{\widehat}

\newcommand{\LB}{\mathcal{L}\mathit{ieb}}

\newcommand{\HoLB}{\mathcal{H}\mathit{olieb}}

\newcommand{\Def}{\mathsf{Def}}

\newcommand{\OOGCdd}{{\mathsf{GC}}^{\mathsf{+or}}_{d,d+1}}
\newcommand{\OGC}{\mathsf{GC}^{\mathsf{or}}}
\newcommand{\OGCd}{\mathsf{GC}^{\mathsf{or}}_{d+1}}

\newcommand{\rOGCd}{\overline{\mathsf{GC}}^{\mathsf{or}}_{d+1}}
\newcommand{\SGCd}{\mathsf{GC}^{\mathsf s}_{d+1}}
\newcommand{\TGCd}{\mathsf{GC}^{\mathsf t}_{d+1}}

\newcommand{\SGC}{\mathsf{GC}^{\mathsf s}}
\newcommand{\TGC}{\mathsf{GC}^{\mathsf t}}
\newcommand{\TGCdd}{\mathsf{GC}^{\mathsf t}_{d,d+1}}
\newcommand{\TTGCdd}{{\mathsf{GC}}^{+\mathsf t}_{d,d+1}}
\newcommand{\rTGCd}{\overline{\mathsf{GC}}_{d+1}^{\mathsf t}}

\newcommand{\STGC}{\mathsf{GC}^{\mathsf{s\cdot t}}}
\newcommand{\STGCd}{\mathsf{GC}^{\mathsf{s\cdot t}}_{d+1}}

\newcommand{\SoTGCd}{\mathsf{GC}^{\mathsf{s+t}}_{d+1}}

\newcommand{\SoTGC}{\mathsf{GC}^{\mathsf{s+t}}}

\newcommand{\GCd}{\mathsf{GC}_d}
\newcommand{\GCD}{\mathsf{GC}_{d+1}}

\newcommand{\dGCd}{\mathsf{dGC}_d}
\newcommand{\dGCD}{\mathsf{dGC}_{d+1}}

\newcommand{\rrdGCD}{\mathsf{d}\overline{\mathsf{GC}}_{d+1}}

\newcommand{\rdGCD}{\mathsf{d}\widetilde{\mathsf{GC}}_{d+1}}

\newcommand{\ddGCdd}{\widehat{\mathsf{dGC}}_{d,d+1}}
\newcommand{\dddGCdd}{{\mathsf{dGC}}^{+}_{d,d+1}}

\newcommand{\st}{\stackrel{\leftrightarrow}}

\newcommand{\OOGCd}{\overline{\mathsf{GC}}^{\mathsf{or}}_{d+1}}

\newcommand{\lGC}{\overline{\GC}_{d+1}^{\mathsf{\lambda}}}

\newcommand{\llGCdd}{{\GC}_{d,d+1}^{\mathsf{\lambda+}}}

\newcommand{\deltabw}{\delta_{\circ\hspace{-0.5mm}-\hspace{-1mm}\bu}}
\newcommand{\lo}{\rightsquigarrow}

\newcommand{\Q}{{\mathbb Q}}
 \newcommand{\Z}{{\mathbb Z}}
 \newcommand{\bS}{{\mathbb S}}

 \newcommand{\K}{{\mathbb K}}

 \newcommand{\ot}{\otimes}

\newcommand{\sX}{{\mathsf X}}
\newcommand{\sY}{{\mathsf Y}}

\newcommand{\GC}{\mathsf{GC}}
\newcommand{\dGC}{\mathsf{dGC}}

\newcommand{\grt}{\fg\fr\ft}
 \newcommand{\Beq}{\begin{equation}}
 \newcommand{\Eeq}{\end{equation}}
 \newcommand{\Beqr}{\begin{eqnarray}}
 \newcommand{\Eeqr}{\end{eqnarray}}
 \newcommand{\Beqrn}{\begin{eqnarray*}}
 \newcommand{\Eeqrn}{\end{eqnarray*}}
 \newcommand{\Ba}{\begin{array}}
 \newcommand{\Ea}{\end{array}}
 \newcommand{\Bi}{\begin{itemize}}
 \newcommand{\Ei}{\end{itemize}}
 \newcommand{\Bc}{\begin{center}}
 \newcommand{\Ec}{\end{center}}

 \newcommand{\fg}{{\mathfrak g}}

\newcommand{\fr}{{\mathfrak r}}

\newcommand{\ft}{{\mathfrak t}}

 \newcommand{\cM}{{\mathcal M}}

 \newcommand{\ga}{\gamma}
 
 \newcommand{\Ga}{\Gamma}

 \newcommand{\la}{\lambda}

 \newcommand{\Img}{{\mathsf I\mathsf m}\, }

 \newcommand{\sip}{\smallskip}
 
 \newcommand{\mip}{\vspace{2.5mm}}

\theoremstyle{plain}
\swapnumbers

\newtheorem{prop-def}[theorem]{Proposition-definition}

\newtheorem{f-theorem}{Formality Theorem}[section]
\newtheorem{main-theorem}{Main~Theorem}[section]
\newtheorem{section-theorem}{Theorem}[section]

\theoremstyle{definition}

\newcommand{\oGCD}{\mathsf{GC}^{\mathsf{T}}_{d}}

\begin{document}

 \sloppy

 \newenvironment{proo}{\begin{trivlist} \item{\sc {Proof.}}}
  {\hfill $\square$ \end{trivlist}}

\long\def\symbolfootnote[#1]#2{\begingroup%
\def\thefootnote{\fnsymbol{footnote}}\footnote[#1]{#2}\endgroup}

 \title{On interrelations between graph complexes}

\author{Sergei A.\ Merkulov}
\address{Sergei~Merkulov:  Department of Mathematics, Luxembourg University,
Maison du Nombre, 6 Avenue de la Fonte,
 L-4364 Esch-sur-Alzette,   Grand Duchy of Luxembourg}
\email{sergei.merkulov@uni.lu}

 \begin{abstract} We study Maxim Kontsevich's graph complex $\GCd$ for any integer $d$  as well as its oriented and targeted versions, and show new short proofs of the theorems due to Thomas Willwacher and Marko \v Zivkovi\' c which establish isomorphisms of their cohomology groups. A new result relating the cohomology of the sourced-targeted graph complex in dimension $d+1$ with the direct sum of {\em two}\, copies of the cohomology group of Maxim Kontsevich's graph complex $\GCd$ in dimension $d$ is obtained. We introduce a new graph complex spanned by {\em purely trivalent}\, graphs and show that its cohomology is isomorphic to $H^\bu(\GCd)$.


\end{abstract}
 \maketitle


{\large
\section{\bf Introduction}
}
\label{sec:introduction}

\subsection{M.Kontsevich's graph complex}
For any integer $d$ there is a remarkable graph complex $\GCd$ which has been introduced  by Maxim Kontsevich
in the context of the theory of deformation quantizations of Poisson structures \cite{Ko1}. Graph complexes $\GCd$ with $d$ of the same parity are isomorphic to each other (up to degree shifts). The cohomology of $\GCd$ is largely unknown except that for $d$ even it has been proven by Thomas Willwacher \cite{Wi1} that $H^{i}(\GC_2)$ vanishes for $i<0$ and for $i=0$ it is isomorphic to the Lie algebra $\grt_1$ of the famous prounipotent Grothendieck-Teichm\"uller group introduced by Vladimir Drinfeld in \cite{Dr},
 $$
 H^0(\GC_2)=\grt_1.
 $$
 The graph complexes $\GCd$ have found many remarkable applications in the deformation quantization theory of Poisson structures and of Lie bialgebras \cite{Ko1,Ko2, Do, MW3, AM}, in the theory of operads and props \cite{Wi1,Wi2,MW1}, in the algebraic topology (see \cite{FTW} and references cited there), in the theory of moduli spaces $\cM_g$ of genus $g$ algebraic curves \cite{CGP, AWZ}, in the Lie theory \cite{AT}, etc. In some of these applications the original Kontsevich's complex $\GCd$ appears not in the original form, but in different (often unexpected and non-obvious) incarnations which have need studied by Thomas Willwacher and Marko \v Zivkovi\' c in \cite{Wi2,Z1,Z2}.

 \subsection{Purely trivalent model of the Kontsevich graph complex}

 We introduce  in \S 3 a new graph complex
  $\oGCD$ generated by {\em purely trivalent}\, graphs (modulo an IHX-type relation) which have edges of two different types; the differential in $\oGCD$ keeps the number of vertices invariant but changes the type of edges. This new complex gives us a trivalent model of $\GC_d$ as we prove the following

  \sip

 {\bf Main Theorem 1}. $H^\bu(\oGCD)=H^\bu(\GC_d)$ {\em for any}\,  $d\in \Z$.

\subsection{Oriented and targeted graph complexes}
The Kontsevich graph complex $\GCd$ is generated, by the very definition,  by graphs with all vertices at least trivalent. It has a directed version $\dGCd$ whose generators are graphs with a {\em fixed direction}\, on edges and  with all vertices at least bivalent and with at least one vertex of valency $\geq 3$. It has been proven by Thomas Willwacher in \cite{Wi1} that a natural morphism of complexes
\Beq\label{1: GC to dGC}
f: \GCd \lon \dGCd
\Eeq
which sends a graph $\Ga$ with no fixed directions on edges into a sum of graphs obtained from $\Ga$ by choosing directions on edges in all possible ways, is a quasi-isomorphism; hence the directed complex $\dGCd$ gives us essentially nothing new at the cohomology level. However the latter complex contains  several interesting (from the point of view of applications) subcomplexes:
\Bi
\item[(i)] let $\OGC_{d}\subset \dGC_d$  be a subcomplex spanned by graphs whose directed edges never form a {\em closed}\, directed path. It is shown in \cite{MW1} that this subcomplex
     --- often called the {\em oriented graph complex}  ---  controls the deformation theory  of
     the properad $\LB_{p,q}$ of Lie bialgebras with $p+q+1=d$, the one in which Lie bracket had degree $1-p$ and Lie cobracket has degree $1-q$ (see \S {\ref{2: subsection on HoLB}} below for details), In particular, the complex $\OGC_3$ it controls the deformation quantization theory of ordinary Lie bialgebras as well as the homotopy theory of their universal quantizations \cite{MW3}.

    \sip

\item[(ii)] Let $\TGC_d$ be a so called {\em targeted graph complex}\, which is by definition  a subcomplex of $\dGC_d$  spanned by graphs having at least one vertex  (called a {\em target}) with no outgoing edges ; this complex is obviously isomorphic
    to the so called {\em sourced}\, subcomplex $\SGC_d \subset \dGC_d$ spanned graphs containing at least one vertex  (called a {\em source})
    with no incoming edges; the latter complex controls the deformation theory of quasi-Lie bialgebras.
\Ei

It has been proven by Thomas Willwacher in \cite{Wi2} that there is an isomorphism
of cohomology groups for any $d\in \Z$,
\Beq\label{1: GC-OGC}
H^\bu(\GCd)=H^\bu(\OGCd).
\Eeq
A second proof of this result has been obtained by  Marko \v Zivkovi\' c in \cite{Z1}; he has also  established isomorphisms of the following cohomology groups
\Beq\label{1: GC-TGC-SGC}
H^\bu(\GCd)=H^\bu(\SGCd)=H^\bu(\TGCd)\ \ \forall d\in \Z.
\Eeq

We show in this paper new proofs of isomorphisms (\ref{1: GC-OGC}) and (\ref{1: GC-TGC-SGC}) using new ``interpolating" graph complexes $\OGC_{d,d+1}$ and $\TGC_{d,d+1}$  spanned by graphs with
{\em two}\, types of vertices --- white ones in the cohomological degree $d$ (as in the case of graphs from $\GCd$) and black ones in the cohomological degree $d+1$ (as in the case of graphs
from $\GC_{d+1}$).

\sip

{\bf Main Theorem 2}. {\em  There exists  graph complexes  $\OGC_{d,d+1}$ and $\TGC_{d,d+1}$ which come equipped with explicit isomorphisms of cohomology groups,}
$$
H^\bu(\GCd)\rightarrow H^\bu(\OGC_{d,d+1}) \leftarrow H^\bu(\OGCd), \ \
H^\bu(\GCd)\rightarrow H^\bu(\TGCdd) \leftarrow H^\bu(\TGCd).
$$

\subsection{Sourced-targeted graph complexes}. Let
$$
\STGCd:=\TGCd \cap \SGCd, \ \ \ \ \SoTGC_d:= \TGCd + \SGCd
$$
be the so called {\em sourced-targeted graph complex}\, and {\em sourced-or-targeted graph complex} respectively. The first complex is shown in \cite{F} to control partially the deformation theory of the {\em wheeled}\, closure  $\HoLB_{p,q}^\circlearrowright$ of the properad $\HoLB_{p,q}$ of strongly homotopy Lie bialgebras with $p+q+1=d+1$. We show in \S 7 of this paper that the associated Mayer-Vietoris short exact sequence
$$
0\lon \STGCd \lon \TGCd \oplus \SGCd \lon \SoTGCd \lon 0
$$
produces a {\em split}\, long exact sequence in cohomology.

\sip

{\bf Main Theorem 3} {\em For any $d\in \Z$ there is a short exact sequence of cohomology groups}
\Beq\label{1: GC-STGC}
0\lon H^{\bu-1}(\SoTGCd) \lon H^\bu(\STGCd)\lon H^\bu(\TGCd) \oplus H^\bu(\SGCd) \lon 0
\Eeq

\mip

This result implies that $H^\bu(\STGCd)$ contains  always  {\em two}\, copies of the cohomology group of the Kontsevich graph complex. In particular, in the most important for applications case $d=3$ it implies isomorphisms of cohomology groups
$$
H^0(\STGC_3)=H^0(\GC_2) \oplus H^0(\GC_2)\simeq \grt_1\oplus \grt_1, \ \ \
H^1(\STGC_3)=H^1(\GC_2) \oplus H^1(\GC_2).
$$
The first of the above two isomorphisms has been established earlier by  Marko \v Zivkovi\' c in \cite{Z3}
using a very short and nice argument based on Thomas Willwacher's estimations of $H^\bu(\GCd)$ (which we remind below in \S 2). Unfortunately, our proof of the above Theorem is not short --- its main part is presented in \S 7 but it uses some other results obtained in \S 4 and  \S 6.

\subsection{Open problems}
There are still open problems left in the theory of the above mentioned graph complexes.
The complexes $\GC_d$, $\TGCd$ and $\STGCd$ are dg {\em Lie algebras}, and it is quite natural to ask (i) a  question
on whether or not the isomorphisms (\ref{1: GC-OGC}), (\ref{1: GC-TGC-SGC}) and the maps in   (\ref{1: GC-STGC}) respect Lie algebra structures, or even (ii) a question on whether or not
there exist zig-zags of quasi-isomorphisms of some {\em dg}\, Lie algebras behind these linear maps. The question (i) has been answered in the affirmative in the case (\ref{1: GC-OGC}) (and hence also in the case (\ref{1: GC-TGC-SGC})) by Thomas Willwacher in \cite{Wi2}, and the proof is quite non-trivial --- it is based on the graph complexes interpretation of the deformation theory of the chain operad of the little disks obtained earlier in the fundamental paper \cite{Wi1} by the same author. The aforementioned open problem (ii) appears to be more complicated
than just proving isomorphisms (\ref{1: GC-OGC}) and (\ref{1: GC-TGC-SGC})  at the level of graded vector spaces.

\subsection{The structure of the paper}  The most important and non-trivial results concerning the structure of the cohomology groups $H^\bu(\GCd)$ are obtained indirectly, through their applications in other mathematical theories \cite{Wi1,Wi2}. This work belongs to a relatively small list of papers  where all new results (except only the new purely trivalent incarnation $\oGCD$ of $\GCd$ discussed in \S 3.3) are obtained from the first principles, i.e.\ using solely the main definitions
of the graph complexes involved. Hence we tried to make the paper self-contained and give in \S 2 an overview of the theory of graph complexes which starts with basic definitions, examples, and an outline of the known results. In \S 3 we remind some basic facts about the deformation theory of (degree shifted) properads of Lie bialgebras and use them to construct a trivalent model $\oGCD$ of the Kontsevich graph complex $\GCd$ and prove Main Theorem 1. Sections 4 is devoted to the study of the directed graph complex $\dGCd$ and its various reduced models which we use in the proofs later.
In \S 5 we introduce auxiliary extensions of $\dGCd$ in terms of graphs with two types of vertices. The latter extensions are used in \S 6 to prove Main Theorem 2. Finally in \S 7 we prove Main Theorem 3.

\subsection{Some notation} We work over a field $\K$ of characteristic zero.
 The set $\{1,2, \ldots, n\}$ is abbreviated to $[n]$;  its group of automorphisms is
denoted by $\bS_n$; the trivial (resp., the sign) one-dimensional representation of
 $\bS_n$ is denoted by $\id_n$ (resp.,  $\sgn_n$). The cardinality of a finite set $S$ is denoted by $\# S$ while its linear span over a
field $\K$ by $\K\left\langle S\right\rangle$.
The top degree skew-symmetric tensor power
of $\K\left\langle S\right\rangle$ is denoted by $\det S$; it is assumed that $\det S$
is a 1-dimensional Euclidean space associated with the unique Euclidean
structure on $\K\left\langle S\right\rangle$ in which the elements of $S$ serve
as an orthonormal basis; in particular, $\det S$ contains precisely two vectors of unit length.
If $V=\oplus_{i\in \Z} V^i$ is a graded vector space, then
$V[k]$ stands for the graded vector space with $V[k]^i:=V^{i+k}$. For $v\in V^i$ we set $|v|:=i$.


\sip

{\bf Acknowledgement}. It is a great pleasure to thank Thomas Willwacher and  Marko \v Zivkovi\'c  for many valuable
discussions and comments. A part of this work has been done while the author visited the IHES (Bures-sur-Yvette), the University of Strasbourg and the ETH (Zurich); I am very grateful to all these institutions for hospitality. I am also grateful to the University of Luxembourg for granting me a sabbatical leave in the summer semester of 2023. Finally I would like to thank an anonymous referee for
providing me with very helpful and insightful remarks and suggestions.

{\large
\section{\bf An overview of the theory of graph complexes}
}
\label{sec:graph complexes}

\subsection{Directed graphs} A directed graph $\Ga$ is a pair of sets $(E(\Ga), V(\Ga))$
called the set of edges and, respectively, the set of vertices of $\Ga$ together
with a map
$$
\Ba{rccc}
h: & E(\Ga) & \lon & V(\Ga)\times V(\Ga)\\
   & e      & \lon & (in_e, out_e)
\Ea
$$
The vertices $in_e$ and $out_e$ are said to be connected by en edge directed from the vertex $in_e$ towards
the vertex  $out_e$. If  $in_e=out_e\in V(\Ga)$, the edge $e$ is called a {\em loop}, while the vertex is called a {\em tadpole}. The number
$$
g(\Ga):= \# E(\Ga) - \#V(\Ga) +1
$$
is called the {\em loop}\, number of $\Ga$. {\em All}\, the differentials in the graph complexes considered in this paper preserve this integer, an important fact  which
 assures  convergence of all spectral sequences used in the arguments (see Appendix G in \cite{Wi1} for full details).

\sip

Every graph $\Ga$ has an obvious geometric interpretation as a 1-dimensional CW complex whose 0-cells are vertices $v\in V(\Ga)$, and 1-cells are edges  $e\in E(\Ga)$; note that each one cell $e$ comes equipped with one of the two possible orientations which is pictorially represented as an arrow (or flow) on that edge.
  Here is an example of a directed graph
 $$
   \Ba{c}\resizebox{14mm}{!}
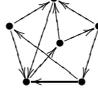
{
\xy
 (1,0)*{\bullet}="0",
(0,8)*{\bullet}="1",
(-8,3)*{\bullet}="5",
(8,3)*{\bullet}="2",
(-5,-7)*{\bullet}="4",
(5,-7)*{\bullet}="3",
\ar @{->} "0";"1" <0pt>
\ar @{->} "0";"2" <0pt>
\ar @{<-} "0";"4" <0pt>
\ar @{->} "1";"4" <0pt>
\ar @{<-} "5";"3" <0pt>
\ar @{<-} "1";"2" <0pt>
\ar @{->} "2";"3" <0pt>
\ar @{->} "3";"4" <0pt>
\ar @{<-} "4";"5" <0pt>
\ar @{->} "5";"1" <0pt>
\endxy}
\Ea
 $$
\vspace{1mm}

A graph is called {\em connected}\, if it is connected as a CW complex.
If a  vertex $v$ of a directed graph has
$m\geq 0$ outgoing edges and $n\geq 0$ incoming edges, then $v$ is called an $(m,n)$-{\em vertex}, and the number $m+n$ is called the {\em valency}\, of $v$. A $(1,1)$-vertex is called {\em passing}.

\subsection{Complexes of directed graphs}

Let $\mathsf{dG}^{\geq 2}_{n,l}$ be the set of connected directed graphs $\Ga$ with $n$ vertices of valency $\geq 2$   (but with no passing vertices, that is bivalent vertices of the form  $
\Ba{c}\resizebox{10mm}{!}{\xy
 (0,1)*{}="a",
(7,1)*{\bu}="b",
 (14,1)*{}="c",
\ar @{->} "a";"b" <0pt>
\ar @{->} "b";"c" <0pt>
\endxy}\Ea$ ) and $l$  edges (but with no loop edges) such that
some bijections $V(\Ga)\rar [n]$ and $E(\Ga)\rar [l]$ are fixed, i.e.\ every edge
and every vertex of $\Ga$ has a fixed numerical label. There is
a natural  action of the group $\bS_n \times  \bS_l$ on the set
$\mathsf{dG}^{\geq 2}_{n,l}$ with $\bS_n$ acting by relabeling the vertices and
 $\bS_l$ acting by relabeling the
edges. 
For a fixed integer $d\in \Z$ consider a graded vector space
\Beq\label{2: dGCd definition}
\dGC_d^{\geq 2}= \prod_{l\geq 2}\prod_{n\geq 2} \K \langle \mathsf{dG}^{\geq 2}_{n,l}\rangle
\ot_{\bS_n\times \bS_l} \left(\sgn_n^{\ot |d|} \ot \sgn_l^{\ot|d-1|}\right) [d(1-n) + l(d-1)]
\Eeq
 This space is spanned by directed graphs with no numerical labels on vertices and
 edges but with a choice (up to a sign) of an {\em orientation}: for $d$ even (resp., odd)
   this is a choice of the ordering of edges $e_1\wedge \ldots \wedge e_l$
   (resp., of vertices $v_1\wedge\ldots \wedge v_n$) up to an even permutation. Each graph $\Ga$ has
precisely two different orientations, $or$ and $-or$; we identify
$(\Ga,or)=-(\Ga,-or)$ and abbreviate the pair $(\Ga,or)$ by $\Ga$.
   The cohomological degree of $\Ga\in \dGC_d^{\geq 2}$ is given by
$$
|\Ga|=d(\# V(\Ga)-1) + (1-d) \#E(\Ga)
$$
or, in terms of the loop number $g:=\# E(\Ga)-\# V(\Ga)+1$ of $\Ga$, one can write
\Beq\label{2: degree of Ga in terms of g}
|\Ga|=-dg + \#E(\Ga)=(1-d)g+\#V(\Ga) -1.
\Eeq

The space $\dGC_d^{\geq 2}$ can be made into a complex \cite{Ko1, Wi1}:
 the differential $\delta$ on $\dGC_d^{\geq 2}$ is defined by its action,
  $$
  \delta\Ga=\sum_{v} \delta_v\Ga,
  $$
  on each vertex
   $$
v= \Ba{c}\resizebox{10mm}{!}{\xy
 (0,0)*{\bullet}="a",
(2.3,4)*{}="1",
(-2.3,4)*{}="2",
(5,1)*{}="3",
(-5,1)*{}="4",
(-3.6,-3.6)*{}="5",
(3.6,-3.6)*{}="6",
(0,-4.5)*{}="7",
\ar @{->} "a";"1" <0pt>
\ar @{->} "a";"2" <0pt>
\ar @{->} "a";"3" <0pt>
\ar @{<-} "a";"4" <0pt>
\ar @{->} "a";"5" <0pt>
\ar @{->} "a";"6" <0pt>
\ar @{<-} "a";"7" <0pt>
%
\endxy}
\Ea \in V(\Ga)
$$
by splitting it into two new vertices $v'$ and $v''$ connected by a
new directed edge $e'$,  and then re-attaching all the half-edges attached
earlier to $v$ to the new vertices
  in all possible ways
\Beq\label{2: delta_v splits vertex into two}
\delta_v:\
\Ba{c}\resizebox{10mm}{!}{\xy
 (0,0)*{\bullet}="a",
(2.3,4)*{}="1",
(-2.3,4)*{}="2",
(5,1)*{}="3",
(-5,1)*{}="4",
(-3.6,-3.6)*{}="5",
(3.6,-3.6)*{}="6",
(0,-4.5)*{}="7",
\ar @{->} "a";"1" <0pt>
\ar @{->} "a";"2" <0pt>
\ar @{->} "a";"3" <0pt>
\ar @{<-} "a";"4" <0pt>
\ar @{->} "a";"5" <0pt>
\ar @{<-} "a";"6" <0pt>
\ar @{<-} "a";"7" <0pt>
%
\endxy}
\Ea
\ \
\lon
\ \
 \sum
\Ba{c}\resizebox{12mm}{!}{\xy
%
%
 (0,-3.3)*{\bullet}="a", (2,-1.7)*{_{v'}},
 (0,3.3)*{\bullet}="b",(2.3,1.7)*{^{v''}},
(-7,-6)*{}="1",
(7,-7)*{}="2",
(-3,-9)*{}="3",
(3,-9)*{}="4",
(5,8)*{}="5",
(-5,8)*{}="6",
(0,9)*{}="7",
\ar @{->} "a";"b" <0pt>
\ar @{->} "a";"1" <0pt>
\ar @{->} "a";"2" <0pt>
\ar @{<-} "a";"3" <0pt>
\ar @{<-} "a";"4" <0pt>
\ar @{->} "b";"5" <0pt>
\ar @{<-} "b";"6" <0pt>
\ar @{->} "b";"7" <0pt>
\endxy}
\Ea
\Eeq
such that no new univalent or passing vertices are created. Here one must be
careful about the signs, i.e.\ about the induced orientations of the summands in
 $\delta_v\Ga$. We adopt the following rule:
\Bi
\item[(i)] for $d$ even the newly created edge $e'$ takes the first position with
respect to the  given ordering of edges of $\Ga$;
\item[(ii)]  for $d$ odd we first choose an ordering $v_1\wedge\ldots\wedge   v_l$
of vertices of $\Ga$ in such a way that $v$ is in the first position;
 then the ordering of vertices in $\delta_v\Ga$ is obtained by
 replacing $v\rar v'\wedge v''$.
\Ei

For future use we abbreviate the action formula for the differential $\delta_v$ as follows
$$
\delta_v: \bu \rightsquigarrow \resizebox{11mm}{!}{
\xy
 (0,1)*{\bullet}="a",
(8,1)*{\bu}="b",
\ar @{->} "a";"b" <0pt>
\endxy}
$$
always assuming that edges attached to $v$ are re-distributed among the two new vertices of the graph
on the r.h.s.\ in all possible ways (which do not contradict its defining conditions).

\sip

Let $\dGCd$ be a subspace of $\dGCd^{\geq 2}$ spanned by graphs having at least one vertex of valency $\geq 3$. The complex $\dGCd^{\geq 2}$ splits as a direct sum of complexes,
$$
\dGCd^{\geq 2}=\dGCd \oplus \dGCd^2
$$
where $\dGCd^2$ is spanned by graphs with all vertices bivalent. The cohomology of the latter complex is not hard to compute (see e.g.\ \cite{Wi1} for details),
$$
H^\bu(\dGCd^2)= \bigoplus_{j\geq 1\  \&\  j\equiv 2d+1 \bmod 4} \K[d-j],
$$
so the most interesting part of the directed graph complex is $\dGCd$.

\sip

The complex $\dGC_d$ has a pre-Lie algebra structure,
$$
\Ba{rccc}
\circ: & \dGC_d \ot \dGCd & \lon & \dGCd\\
     & \Ga_1\ot \Ga_2 & \lon & \Ga_1\circ \Ga_2:=\sum_{v\in V(\Ga_1)} \Ga_1 \circ_v \Ga_2,
\Ea
$$
where $\Ga_1 \circ_v \Ga_2$ is obtained from $\Ga_1$ by substituting into its vertex $v$ the graph  $\Ga_2$
and taking the sum over all possible re-attachments of edges connected earlier to $v$ to vertices of $\Ga_2$ in all possible ways. Using this operation the differential $\delta$ can be described as the composition
$$
\delta\Ga=(-1)^{|\Ga|} \Ga\circ  \left( \resizebox{11mm}{!}{
\xy
 (0,1)*{\bullet}="a",
(8,1)*{\bu}="b",
\ar @{->} "a";"b" <0pt>
\endxy}\right) \ \text{modulo graphs with univalent vertices}.
$$
The operation $\circ$  is compatible with the differential so that $\dGC_d$ is a dg Lie algebra with the Lie bracket given by
$$
[\Ga_1,\Ga_2]:= \Ga_1\circ \Ga_2 - (-1)^{|\Ga_1||\Ga_2|} \Ga_2\circ \Ga_1.
$$

\subsection{Complexes of undirected graphs}
Let $\mathsf{dG}^{ 2}_{n,l}$ be  the set of directed connected graphs with every vertex at least bivalent and with no loops (passing vertices are now allowed).
If we split an edge of a graph $\Ga\in \mathsf{G}^{2}_{n,l}$ into two half-edges
$(h_1, h_2)$ (e.g. by cutting the edge $e$ in the middle using its geometric incarnation), then a direction of $e$ can be understood as a total ordering of the set $H(e):=(h_1,h_2)$, that is a choice of one of the two possible unital vectors in the 1-dimensional vector space
$\det H(e)$. We use this observation below.

\sip

There is
a natural action of  the group
$\bS_n\times  (\bS_l \ltimes (\bS_2)^l)$ on the set $\mathsf{dG}^{\geq 2}_{n,l}$ with
$(\bS_2)^l$ acting by reversing the directions of the edges. For a fixed integer
$d\in \Z$  consider a graded vector space,
$$
\GC_d^{\geq 2}:= \prod_{l\geq 3}\prod_{n\geq 2} \K \langle \mathsf{dG}^2_{n,l}\rangle
\ot_{\bS_n\times  (\bS_l \ltimes
 (\bS_2)^l)} \left(\sgn_n^{\ot |d|}
\ot \sgn_l^{|d-1|}\ot \sgn_2^{\ot |dl|}\right) [d(1-n) + l(d-1)],
$$
which can be understood as generated by at least bivalent graphs $\Ga$ with no
labels on vertices and edges but equipped instead with an orientation $or$ defined
as a length 1 vector $or$ of the following Euclidean one-dimensional space,
$$
or \in \left\{\Ba{cc} \det E(\Ga)  & \text{\sf if $d$ is even},\\
\det V(\Ga)
\bigotimes_{e\in E(\Ga)} \det H(e) & \text{\sf if $d$ is odd}.
\Ea
\right.
$$
For $d$ even the edges of graphs can be understood as undirected, and their orientation
becomes just the ordering of edges up to an even permutation.
 For $d$ odd
the orientation $or$ of $\Ga$ can be equivalently understood as an ordering
 of its vertices (up to a sign), and a
choice of the direction on each edge $e$, up to a flip and multiplication by $(-1)^d$.
We make edges of graphs from $\GCd^{\geq 2}$ {\em dotted}\, in our pictures in order not to confuse such edges with solid ones whose directions are strictly fixed; moreover we often (but not always) omit showing arrows
on dotted edges as their direction can be flipped,
$\xy
 (0,0)*{\bullet}="a",
(7,0)*{\bu}="b",
\ar @{.>} "a";"b" <0pt>
\endxy = (-1)^{d} \xy
 (0,0)*{\bullet}="a",
(7,0)*{\bu}="b",
\ar @{<.} "a";"b" <0pt>
\endxy$.

\sip

The differential $\delta$ on $\GC_d^{\geq 2}$ is given by its action on each vertex
$
v= \Ba{c}\resizebox{7mm}{!}{\xy
 (0,0)*{\bullet}="a",
(2.3,4)*{}="1",
(-2.3,4)*{}="2",
(5,1)*{}="3",
(-5,1)*{}="4",
(-3.6,-3.6)*{}="5",
(3.6,-3.6)*{}="6",
(0,-4.5)*{}="7",
\ar @{.} "a";"1" <0pt>
\ar @{.} "a";"2" <0pt>
\ar @{.} "a";"3" <0pt>
\ar @{.} "a";"4" <0pt>
\ar @{.} "a";"5" <0pt>
\ar @{.} "a";"6" <0pt>
\ar @{.} "a";"7" <0pt>
\endxy}
\Ea
 $
 of a graph $\Ga\in \GC_d^{\geq 2}$
by splitting $v$ into two new vertices connected by an edge,
$$
 \delta\Ga=\sum_v \delta_v\Ga, \ \ \ \
\delta_v:
 \Ba{c}\resizebox{9mm}{!}{\xy
 (0,0)*{\bullet}="a",
(2.3,4)*{}="1",
(-2.3,4)*{}="2",
(5,1)*{}="3",
(-5,1)*{}="4",
(-3.6,-3.6)*{}="5",
(3.6,-3.6)*{}="6",
(0,-4.5)*{}="7",
\ar @{.} "a";"1" <0pt>
\ar @{.} "a";"2" <0pt>
\ar @{.} "a";"3" <0pt>
\ar @{.} "a";"4" <0pt>
\ar @{.} "a";"5" <0pt>
\ar @{.} "a";"6" <0pt>
\ar @{.} "a";"7" <0pt>
\endxy}
\Ea
\ \lon \ \sum
\Ba{c}\resizebox{10mm}{!}{\xy
%
%
 (0,-2.3)*{\bullet}="a",
 (0,2.3)*{\bullet}="b",
(-7,-5)*{}="1",
(7,-6)*{}="2",
(-3,-8)*{}="3",
(3,-8)*{}="4",
(5,7)*{}="5",
(-5,7)*{}="6",
(0,8)*{}="7",
\ar @{.} "a";"b" <0pt>
\ar @{.} "a";"1" <0pt>
\ar @{.} "a";"2" <0pt>
\ar @{.} "a";"3" <0pt>
\ar @{.} "a";"4" <0pt>
\ar @{.} "b";"5" <0pt>
\ar @{.} "b";"6" <0pt>
\ar @{.} "b";"7" <0pt>
\endxy}
\Ea \ \ \ (\text{or, shortly}, \ \delta_v: \bu \rightsquigarrow \resizebox{11mm}{!}{
\xy
 (0,1)*{\bullet}="a",
(7,1)*{\bu}="b",
\ar @{.} "a";"b" <0pt>
\endxy})
$$
and then re-attaching the edges attached earlier to $v$ to the new vertices
in all possible ways
such that no new vertices of valency $1$ are created.
Orientations of each summand in $\delta_v\Ga$ are defined in the full analogy to
the case $\dGC_d^{\geq 2}$ discussed above.

\sip

The complex $\GCd^{\geq 2}$ contains two subcomplexes as direct summands,
$\GCd^2$ generated by graphs with all vertices bivalent and $\GCd$ with all vertices at least trivalent. One has \cite{Wi1}
$$
H^\bu(\GCd^{\geq 2})= H^\bu(\GCd^2) \oplus H^\bu(\GCd)
$$
and
$$
H^\bu(\GCd^2)=H^\bu(\dGC_d^2)= \bigoplus_{j\geq 1\atop j\equiv 2d+1 \mod 4} \K[d-j],
$$
where the summand $\K[d-j]$ is generated by the polytope-like graph with $j$ vertices,
that is, a connected graph with $j$ bivalent vertices, e.g.\
$\Ba{c}\resizebox{6mm}{!}{\xy
 (0,0)*{\bullet}="a",
(6,0)*{\bu}="b",
(3,5)*{\bu}="c",
\ar @{.} "a";"b" <0pt>
\ar @{.} "a";"c" <0pt>
\ar @{.} "c";"b" <0pt>
\endxy}\Ea$ for $j=3$ and  $d$ odd (for $d$ even this particular graph vanishes identically as it admits an automorphism reversing its orientation). The complex $\GCd$ is in fact a dg Lie algebra
with the Lie bracket defined in a full analogy to the case $\dGC_d$ discussed above.

\sip

The complex $\GCd$ for $d=2$ was introduced and studied by Maxim Kontsevich in \cite{Ko1} in the context of the theory of the deformation quantization of Poisson structures.

\sip


 It is a very hard problem to compute the cohomology groups
$H^\bu(\GC_d)$.
 A remarkable  result \cite{Wi1} by T.\ Willwacher identifies the
zero-th cohomology group $H^0(\GC_2)$ with the Lie algebra $\grt_1$ of the
Grothendieck-Teichm\"uller group $GRT_1$ introduced by  V.\ Drinfeld  in \cite{Dr}.
This group occurs mysteriously in several very different areas of pure
mathematics and mathematical physics. Here is a simplest possible example of a non-trivial cohomology
class in $\GC_2$,
$$
 \Ba{c}\resizebox{12mm}{!}{
\xy
 (0,0)*{\bullet}="a",
(0,8)*{\bullet}="b",
(-7.5,-4.5)*{\bullet}="c",
(7.5,-4.5)*{\bullet}="d",
\ar @{.} "a";"b" <0pt>
\ar @{.} "a";"c" <0pt>
\ar @{.} "b";"c" <0pt>
\ar @{.} "d";"c" <0pt>
\ar @{.} "b";"d" <0pt>
\ar @{.} "d";"a" <0pt>
\endxy}
\Ea\in H^0(\GC_2),
$$
first detected in \cite{Ko1}.

\sip

 There is a {\em quasi-isomorphism}\,  of dg Lie algebras
\cite{Wi1}
\Beq\label{2: GC to dGC}
f: \GCd \lon \dGCd
\Eeq
which sends a graph from $\GC_d$ into a sum of directed graphs by assigning to each edge
a fixed direction in both possible ways, $\resizebox{9mm}{!}{
\xy
 (0,1)*{\bullet}="a",
(8,1)*{\bu}="b",
\ar @{.>} "a";"b" <0pt>
\endxy} \rightsquigarrow \resizebox{9mm}{!}{
\xy
 (0,1)*{\bullet}="a",
(8,1)*{\bu}="b",
\ar @{->} "a";"b" <0pt>
\endxy} + (-1)^d
\resizebox{9mm}{!}{
\xy
 (0,1)*{\bullet}="a",
(8,1)*{\bu}="b",
\ar @{<-} "a";"b" <0pt>
\endxy}
$.

\sip

One can assume without (much) loss of generality that $\GC_d$ is spanned by graphs with
no multiple edges between the same pair of vertices \cite{WZ}. Complexes $\GC_d$
with the parameters $d$ of the same parity are isomorphic to each other
(up to degree shifts). Hence there are essentially only two really different graph
complexes, $\GC_d$ and $\GC_{d+1},$ whatever $d$ is. In applications, however, the
particular values of the parameter $d$ play sometimes an important role; e.g. in applications of graph complexes in mathematical physics and the deformation quantization theories
the parameter $d$ stands often for the dimension of the space in which a quantum
field theory or appropriate configuration spaces are studied.

\sip

It is worth noting that the differential $\delta$ on the complexes $\dGC_d$ and $\GC_d$ (and hence in their subcomplexes described below) does not change the loop numbers $g$ of graphs so the complexes decompose into products over the parameter $g$
$$
\GCd=\prod_{g\geq 2} \GC_{g,d}\ , \ \ \ \dGC_d=\prod_{g\geq 2} \dGC_{g,d},
$$
where $\GC_{g,d}$ and $\dGC_{g,d}$ are subcomplexes spanned by graphs $\Ga$ with loop number $g$.
 It is worth noting that
\Beq\label{2: V(Ga) leq 2g-2}
\#V(\Ga)\leq 2g-2\ \ \ \text{\em for any connected graph $\Ga$ with all vertices at least trivalent}.
\Eeq
Indeed, in this case $2\#E(\Ga)\geq 3\# V(\Ga)$ implying the above estimation.

\sip

Thomas Willwacher has proven the following estimations
on the non-triviality of cohomology groups of the Kontsevich graph complex:

\subsubsection{\bf Theorem \cite{Wi1,Wi3}}\label{2: lemma on vanishing H(GC)}
 {\em  For $g \geq 2$ the graph cohomology $H^k(\GC_{d, g})$ vanishes outside of the range
$$
 (2 - d)g - 1 \leq  k  \leq  (3 - d)g - 3
$$
For $g = 3$ the lower bound can be improved  to $(2 - d)g$.
}

\subsection{Sourced and targeted graph complexes} The quasi-isomorphism  (\ref{2: GC to dGC}) says that at the cohomology level the directed graph complex $\dGC_d$ gives us nothing new comparing to the original Maxim Kontsevich graph complex $\GCd$. However the complex $\dGC_d$ is useful as it contains several interesting subcomplexes  which we discuss in the rest of this section.

\sip

Let $\Ga$ be a graph in $\dGC_d$. A vertex of $\Ga$ is called a {\em target}\, (resp. a {\em source}) if it has no outgoing (resp.\ ingoing) edges. The linear subspace $\GC_d^{\mathsf t}$ (resp., $\GC_d^{\mathsf s}$) of $\dGC_d$
spanned by directed graphs  with at least one target (resp., source) is a dg Lie subalgebra
called the {\em targeted}\, (resp., {\em sourced}\,) graph complex. It has been proven in \cite{Z2} that the cohomology isomorphisms (\ref{1: GC-TGC-SGC}) hold true. We show a new (relatively short) proof of this isomorphism in \S 6 below.

\sip

The directed graph complex $\dGC_d$ admits an involution $\iota$ given by reversing arrows
on its edges; more precisely, the involution $\iota$ is given by the formula \cite{MZ}
\Beq\label{3: involution on dGG_d}
\Ba{rccc}
\imath: & \dGC_d & \lon & \dGC_d\\
        & \Ga & \lon  &
\left\{\Ba{cc} (-1)^{\#E(\Ga)+\#V(\Ga)+1}\st{\Ga}  & \text{for $d$ even} \\
(-1)^{\#V(\Ga) +1} \st{\Ga} & \text{for $d$ odd}
\Ea
\right.
\Ea
\Eeq
 where $\st{\Ga}$ is the directed graph obtained from a graph $\Ga$ by reversing simultaneously directions of its edges while keeping its orientation $or$ unchanged (which makes sense as $V(\st{\Ga})=V(\Ga)$ and $E(\st{\Ga})=E(\Ga)$). The map $\iota$
it commutes with the differential and essentially identifies complexes $\GC_d^{\mathsf t}$ and $\GC_d^{\mathsf s}$.

\sip

The complex $\STGCd=\SGCd\cap \TGCd$ is of interest because it controls \cite{F} the deformation theory of the {\em wheeled}\, closure  $\HoLB_{p,q}^\circlearrowright$ of the properad $\HoLB_{p,q}$ of strongly homotopy Lie bialgebras with $p+q=d-1$ (we remind definitions in \S 3 below). It fits the following short exact sequence
of complexes,
\Beq\label{2: seq STGC SGC sum TGC to SoTGC}
0\lon \STGCd \lon \SGCd \oplus \TGCd \lon \SoTGCd \lon 0
\Eeq
where $\SoTGCd$ is the subcomplex of $\dGCD$ spanned by graphs with at least one source {\em or}\, with at least one target, i.e.\ $\SoTGCd:=\SGCd + \TGCd\subset \dGCD$.
Denote by $\GC_{g,d+1}^{\mathsf{s\cdot t}}$ and $\GC_{g,d+1}^{\mathsf{s+t}}$
be subcomplexes of, respectively, $\STGCd$ and $\SoTGCd$ spanned by graphs with loop number $g$.
Using estimations {\ref{2: lemma on vanishing H(GC)}} and the isomorphism (\ref{1: GC to dGC}) Marko \v Zivkovi\' c  has proven  \cite{Z3} that for $g\geq 3$ and for degrees $k \leq (d - 2)g$ one has the following isomorphism of cohomology groups,
$$
H^k(\GC_{g,d+1}^{\mathsf{s\cdot t}})=H^k(\GCd) \oplus H^k(\GCd),
$$
while for $g\geq 3$ and  for degrees $k > (d - 2)g$ it holds that
$$
H^k(\GC_{g,d+1}^{\mathsf{s\cdot t}})=H^k(\GC^{\mathsf{s+t}}_{g,d+1})
$$

\sip

This result says that in the important for applications case $d=3$ one has,
\Beq\label{2: H^0(STGC) as a sum of grt}
H^0(\STGC_3)=H^0(\GC_2) \oplus H^0(\GC_2)\simeq \grt_1\oplus \grt_1.
\Eeq
We strengthen this nice observation in \S 7 below by proving that there is a short exact sequence of cohomology groups (\ref{1: GC-STGC}) which in turn implies (\ref{2: H^0(STGC) as a sum of grt}) and also an isomorphism of the first
cohomology groups,
$$
H^1(\STGC_3)=H^1(\GC_2) \oplus H^1(\GC_2).
$$



\subsection{Oriented  graph complexes}
Let $\GC^{\mathsf{or,\geq 2}}_d\subset \dGC_d^{\mathsf{\geq 2}}$ be the sub-complex (in fact a dg Lie subalgebra), spanned by directed graphs  with
with no closed paths of directed edges (i.e.\ no so called wheels), e.g.
$$
  \Ba{c}\resizebox{12mm}{!}{
\xy
 (0,0)*{\bullet}="a",
(0,8)*{\bullet}="b",
(-7.5,-4.5)*{\bullet}="c",
(7.5,-4.5)*{\bullet}="d",
\ar @{->} "a";"b" <0pt>
\ar @{<-} "a";"c" <0pt>
\ar @{<-} "b";"c" <0pt>
\ar @{<-} "b";"d" <0pt>
\ar @{->} "d";"a" <0pt>
\endxy}
\Ea \in \GC^{\mathsf{ or,\geq 2}}_d, \ \ \ \
  \Ba{c}\resizebox{12mm}{!}{
\xy
 (0,0)*{\bullet}="a",
(0,8)*{\bullet}="b",
(-7.5,-4.5)*{\bullet}="c",
(7.5,-4.5)*{\bullet}="d",
\ar @{<-} "a";"b" <0pt>
\ar @{->} "a";"c" <0pt>
\ar @{<-} "b";"c" <0pt>
\ar @{<-} "b";"d" <0pt>
\ar @{->} "d";"a" <0pt>
\endxy}
\Ea \notin \GC^{\mathsf{or,\geq 2}}_d.
$$
 It is often called an {\em oriented}\, graph complex in the sense of the fixed
 flow on all edges of $\Ga\in \GC^{or}_d$ directed from, say, bottom of the graph
 to the top. 
It was proven in \cite{Wi2} that
$$
H^\bu(\GC_d^{\mathsf{\geq 2}})=H^\bu(\GC^{\mathsf{or, \geq 2}}_{d+1}) \ \ \ \forall d\in \Z.
$$
Moreover, this isomorphism has been established by Thomas Willwacher in \cite{Wi2} at the level of Lie algebras, not just at the level of graded vector spaces. Let us denote
by $\GC^{\mathsf{or}}_d\subset \GC^{\mathsf{or,\geq 2}}_d$ the subcomplex spanned by graphs with at least one vertex of valency $\geq 3$. Then one has
\Beq\label{2: H(GCd) H(OGCd+1)}
H^\bu(\GC_d)=H^\bu(\GC^{\mathsf{or}}_{d+1}) \ \ \ \forall d\in \Z.
\Eeq

\sip

Marko \v Zivkovi\'c has found a different proof of the isomorphism (\ref{2: H(GCd) H(OGCd+1)}) at the level of graded vector spaces by giving an {\em explicit}\, (and beautiful) morphism of dual graph complexes (not respecting, however, the Lie co-algebra structures),
$$
(\GCd)^* \lon (\OOGCd)^*
$$
where  $\OOGCd$ is a reduced  version  of $\OGCd$ which has been introduced
by Marko  \v Zivkovi\' c in \cite{Z1} and
whose definition  we remind next.

\sip

 By definition \cite{Z1}, the complex $\OOGCd$
is generated by  graphs with {\em all}\, vertices at
 least trivalent, and with two types of edges, solid edges
$
 \xy
 (0,0)*{\bullet}="a",
(7,0)*{\bu}="b",
\ar @{->} "a";"b" <0pt>
\endxy$     and dotted ones
 $\xy
 (0,0)*{\bullet}="a",
(7,0)*{\bu}="b",
\ar @{.>} "a";"b" <0pt>
\endxy.
$
 Solid edges are assigned degree $-d$ and have a fixed direction; dotted edges are assigned degree $1-d$ and their directions
 can be flipped according to the rule
\Beq\label{2: symmetry of dotted edges}
 \xy
 (0,0)*{\bullet}="a",
(7,0)*{\bu}="b",
\ar @{.>} "a";"b" <0pt>
\endxy = (-1)^{d} \xy
 (0,0)*{\bullet}="a",
(7,0)*{\bu}="b",
\ar @{<.} "a";"b" <0pt>
\endxy.
\Eeq
Often we do not show directions of dotted edges in our pictures;
moreover, for $d$ even dotted edges should be viewed as undirected ones in the full
agreement with this abuse of notation.  The most important condition on every
graph $\Ga$ in $\OOGCd$ is that $\Ga$ is {\em not}\, allowed to have closed
paths of directed solid edges; the dotted edges are assumed to be
``non-passable"\, in any direction, i.e.\ no restrictions are imposed on
these edges. For example,

$$
 \Ba{c}\resizebox{15mm}{!}
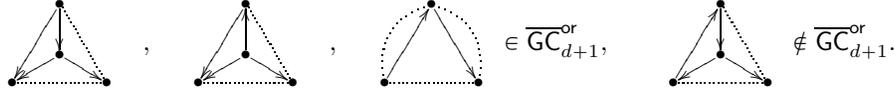
{
\xy
 (0,0)*{\bullet}="a",
(0,8)*{\bullet}="b",
(-7.5,-4.5)*{\bullet}="c",
(7.5,-4.5)*{\bullet}="d",
\ar @{<-} "a";"b" <0pt>
\ar @{->} "a";"c" <0pt>
\ar @{->} "b";"c" <0pt>
\ar @{.} "d";"c" <0pt>
\ar @{.} "b";"d" <0pt>
\ar @{<-} "d";"a" <0pt>
\endxy}
\Ea
\ \ , \ \
 \Ba{c}\resizebox{15mm}{!}{
\xy
 (0,0)*{\bullet}="a",
(0,8)*{\bullet}="b",
(-7.5,-4.5)*{\bullet}="c",
(7.5,-4.5)*{\bullet}="d",
\ar @{->} "a";"b" <0pt>
\ar @{->} "a";"c" <0pt>
\ar @{->} "b";"c" <0pt>
\ar @{.} "d";"c" <0pt>
\ar @{.} "b";"d" <0pt>
\ar @{<-} "d";"a" <0pt>
\endxy}
\Ea \ \ , \ \
  \Ba{c}\resizebox{15mm}{!}{
\xy
 (0,8)*{\bullet}="a",
(-7.5,-4.5)*{\bullet}="c",
(7.5,-4.5)*{\bullet}="d",
\ar@/^-0.8pc/@{.}"a";"c" <0pt>
\ar@/^0.8pc/@{.}"a";"d" <0pt>
\ar @{<-} "a";"c" <0pt>
\ar @{.} "d";"c" <0pt>
\ar @{->} "a";"d" <0pt>
\endxy}
\Ea\in \OOGCd, \ \ \ \
 \Ba{c}\resizebox{15mm}{!}{
\xy
 (0,0)*{\bullet}="a",
(0,8)*{\bullet}="b",
(-7.5,-4.5)*{\bullet}="c",
(7.5,-4.5)*{\bullet}="d",
\ar @{<-} "a";"b" <0pt>
\ar @{->} "a";"c" <0pt>
\ar @{<-} "b";"c" <0pt>
\ar @{.} "d";"c" <0pt>
\ar @{.} "b";"d" <0pt>
\ar @{<-} "d";"a" <0pt>
\endxy}
\Ea   \notin \OOGCd.
 $$
 An orientation $or$ on each $\Ga\in \OOGCd$ is defined
as a length 1 vector of the following Euclidean one-dimensional space
$$
or \in \left\{\Ba{ll} \det E_{sol}(\Ga)\bigotimes_{e\in E_{dot}(\Ga)} \det H(e)   &
\text{ if $d$ is odd},\\
 \det V(\Ga)\ot  \det E_{dot}(\Ga)
\bigotimes_{e\in E_{sol}(\Ga)} \det H(e)   & \text{if $d$ is even}.
\Ea
\right.
$$

The differential on $\OOGCd$ consists of two parts,
\Beq\label{2: delta on hatOGCd}
\delta \Ga =\underbrace{\sum_{v\in V(\Ga)} \delta_{v} \Ga}_{\delta_{\bu}\Ga} +
\underbrace{\sum_{e\in E_{sol}(\Ga)} \delta_e \Ga}_{\delta' \Ga},
\Eeq
where the first part acts on vertices of $\Ga$ by splitting them in a full analogy
to the differential on $\dGC_d$,
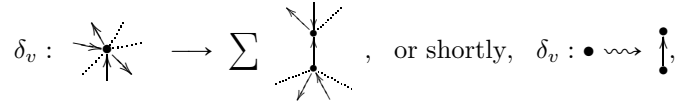
$$
\delta_{v}:
 \Ba{c}\resizebox{9mm}{!}{\xy
 (0,0)*{\bullet}="a",
(2.3,4)*{}="1",
(-2.3,4)*{}="2",
(5,1)*{}="3",
(-5,1)*{}="4",
(-3.6,-3.6)*{}="5",
(3.6,-3.6)*{}="6",
(0,-4.5)*{}="7",
\ar @{.} "a";"1" <0pt>
\ar @{->} "a";"2" <0pt>
\ar @{.} "a";"3" <0pt>
\ar @{<-} "a";"4" <0pt>
\ar @{.} "a";"5" <0pt>
\ar @{->} "a";"6" <0pt>
\ar @{<-} "a";"7" <0pt>
\endxy}
\Ea
\ \lon \ \sum
\Ba{c}\resizebox{10mm}{!}{\xy
%
%
 (0,-3.3)*{\bullet}="a", 
 (0,3.3)*{\bullet}="b",
(-7,-6)*{}="1",
(7,-7)*{}="2",
(-3,-9)*{}="3",
(3,-9)*{}="4",
(5,8)*{}="5",
(-5,8)*{}="6",
(0,9)*{}="7",
\ar @{->} "a";"b" <0pt>
\ar @{.} "a";"1" <0pt>
\ar @{.} "a";"2" <0pt>
\ar @{->} "a";"3" <0pt>
\ar @{<-} "a";"4" <0pt>
\ar @{.} "b";"5" <0pt>
\ar @{->} "b";"6" <0pt>
\ar @{<-} "b";"7" <0pt>
\endxy}
\Ea,\ \  \text{or shortly},\ \
\delta_v: \bu \rightsquigarrow \resizebox{2.5mm}{!}{
\xy
 (0,-2)*{\bullet}="a",
(0,4)*{\bu}="b",
\ar @{->} "a";"b" <0pt>
\endxy},
$$
while the second part  $\delta'$ acts only on solid edges of $\Ga$ by making them dotted,
$$
\delta_e:   \Ba{c}\resizebox{11mm}{!}{
\xy
(0,2)*{^{v_1}},
(8,2)*{^{v_2}},
 (0,0)*{\bullet}="a",
(8,0)*{\bu}="b",
\ar @{->} "a";"b" <0pt>
\endxy}\Ea
\lon
 \Ba{c}\resizebox{11mm}{!}{
\xy
(0,2)*{^{v_1}},
(8,2)*{^{v_2}},
 (0,0)*{\bullet}="a",
(8,0)*{\bu}="b",
\ar @{.>} "a";"b" <0pt>
\endxy}\Ea.
$$
It has been proven in \cite{Z1} that there is a monomorphism of complexes,
$$
j: \OOGCd \lon \GC^{or}_{d+1}
$$
which preserves vertices and solid edges of the generating graphs
$\Ga\in \OOGCd$ while sends every dotted edge of $\Ga$ into the following linear combination
\Beq\label{2: map j from dooted edge to solid comb}
j:\
\Ba{c}\resizebox{10mm}{!}{\xy
(0,2)*{^{v_1}},
(8,2)*{^{v_2}},
 (0,0)*{\bullet}="a",
(8,0)*{\bu}="b",
\ar @{.>} "a";"b" <0pt>
\endxy}\Ea
\lon
\frac{1}{2}\left(
\Ba{c}\resizebox{15mm}{!}{  \xy
(0,2)*{^{v_1}},
(7,5)*{^{v}},
(14,2)*{^{v_2}},
 (0,0)*{\bullet}="a",
 (7,3)*{\bu}="c",
(14,0)*{\bu}="b",
\ar @{->} "a";"c" <0pt>
\ar @{->} "b";"c" <0pt>
\endxy}\Ea
-
\Ba{c}\resizebox{15mm}{!}{
\xy
(0,2)*{^{v_1}},
(7,-1)*{^{v}},
(14,2)*{^{v_2}},
 (0,0)*{\bullet}="a",
 (7,-3)*{\bu}="c",
(14,0)*{\bu}="b",
\ar @{<-} "a";"c" <0pt>
\ar @{<-} "b";"c" <0pt>
\endxy}\Ea
\right)
\Eeq
where the new solid edges on the r.h.s\ are ordered from left to right, while the new vertex
$v$ is assumed to be adjoined into the given ordering of vertices of $\Ga$ as the last one $v_1\wedge v_2\wedge v$. We shall use the graph complex $\OOGCd$ in \S 6 below. In the next section \S 3 we discuss much smaller versions of $\OOGCd$  whose cohomology groups  are still isomorphic to $H^\bu(\GCd)$. The results of that section are not used in the rest of the paper.

\mip


{\large
\section{\bf Reduced versions of the oriented graph complex}
}

\subsection{Graphs complexes and the deformation theory of Lie bialgebras }\label{2: subsection on HoLB} In the arguments below we use a few basic facts about the properad of Lie bialgebras. They are well-known but let us remind them briefly in this subsection for completeness (see e.g. \cite{MW1} and references cited there for full details).
The properad of degree shifted Lie bialgebras $\LB_{p,q}$ is generated by two operations of degrees $1-p$ and $1-q$ respectively,
$$
\Ba{c}\begin{xy}
 <0mm,0.66mm>*{};<0mm,3mm>*{}**@{-},
 <0.39mm,-0.39mm>*{};<2.2mm,-2.2mm>*{}**@{-},
 <-0.35mm,-0.35mm>*{};<-2.2mm,-2.2mm>*{}**@{-},
 <0mm,0mm>*{\bu};<0mm,0mm>*{}**@{},
   <0.39mm,-0.39mm>*{};<2.9mm,-4mm>*{^{_2}}**@{},
   <-0.35mm,-0.35mm>*{};<-2.8mm,-4mm>*{^{_1}}**@{},
\end{xy}\Ea
=(-1)^{q}
\Ba{c}\begin{xy}
 <0mm,0.66mm>*{};<0mm,3mm>*{}**@{-},
 <0.39mm,-0.39mm>*{};<2.2mm,-2.2mm>*{}**@{-},
 <-0.35mm,-0.35mm>*{};<-2.2mm,-2.2mm>*{}**@{-},
 <0mm,0mm>*{\bu};<0mm,0mm>*{}**@{},
   <0.39mm,-0.39mm>*{};<2.9mm,-4mm>*{^{_1}}**@{},
   <-0.35mm,-0.35mm>*{};<-2.8mm,-4mm>*{^{_2}}**@{},
\end{xy}\Ea
\ \ ,\ \
\Ba{c}\begin{xy}
 <0mm,-0.55mm>*{};<0mm,-2.5mm>*{}**@{-},
 <0.5mm,0.5mm>*{};<2.2mm,2.2mm>*{}**@{-},
 <-0.48mm,0.48mm>*{};<-2.2mm,2.2mm>*{}**@{-},
 <0mm,0mm>*{\bu};<0mm,0mm>*{}**@{},
 <0.5mm,0.5mm>*{};<2.7mm,2.8mm>*{^{_2}}**@{},
 <-0.48mm,0.48mm>*{};<-2.7mm,2.8mm>*{^{_1}}**@{},
 \end{xy}\Ea
=(-1)^{p}
\Ba{c}\begin{xy}
 <0mm,-0.55mm>*{};<0mm,-2.5mm>*{}**@{-},
 <0.5mm,0.5mm>*{};<2.2mm,2.2mm>*{}**@{-},
 <-0.48mm,0.48mm>*{};<-2.2mm,2.2mm>*{}**@{-},
 <0mm,0mm>*{\bu};<0mm,0mm>*{}**@{},
 <0.5mm,0.5mm>*{};<2.7mm,2.8mm>*{^{_1}}**@{},
 <-0.48mm,0.48mm>*{};<-2.7mm,2.8mm>*{^{_2}}**@{},
 \end{xy}\Ea,
$$
which represent Lie bracket and, respectively Lie cobracket (i.e.\ the direction flow in these and similar pictures is assumed to go from the bottom to the top). These generators are assumed to satisfy the following relations
$$
\Ba{c}\resizebox{8mm}{!}{
\begin{xy}
 <0mm,0mm>*{\bu};<0mm,0mm>*{}**@{},
 <0mm,-0.49mm>*{};<0mm,-3.0mm>*{}**@{-},
 <0.49mm,0.49mm>*{};<1.9mm,1.9mm>*{}**@{-},
 <-0.5mm,0.5mm>*{};<-1.9mm,1.9mm>*{}**@{-},
 <-2.3mm,2.3mm>*{\bu};<-2.3mm,2.3mm>*{}**@{},
 <-1.8mm,2.8mm>*{};<0mm,4.9mm>*{}**@{-},
 <-2.8mm,2.9mm>*{};<-4.6mm,4.9mm>*{}**@{-},
   <0.49mm,0.49mm>*{};<2.7mm,2.3mm>*{^3}**@{},
   <-1.8mm,2.8mm>*{};<0.4mm,5.3mm>*{^2}**@{},
   <-2.8mm,2.9mm>*{};<-5.1mm,5.3mm>*{^1}**@{},
 \end{xy}}\Ea
 +
\Ba{c}\resizebox{8mm}{!}{\begin{xy}
 <0mm,0mm>*{\bu};<0mm,0mm>*{}**@{},
 <0mm,-0.49mm>*{};<0mm,-3.0mm>*{}**@{-},
 <0.49mm,0.49mm>*{};<1.9mm,1.9mm>*{}**@{-},
 <-0.5mm,0.5mm>*{};<-1.9mm,1.9mm>*{}**@{-},
 <-2.3mm,2.3mm>*{\bu};<-2.3mm,2.3mm>*{}**@{},
 <-1.8mm,2.8mm>*{};<0mm,4.9mm>*{}**@{-},
 <-2.8mm,2.9mm>*{};<-4.6mm,4.9mm>*{}**@{-},
   <0.49mm,0.49mm>*{};<2.7mm,2.3mm>*{^2}**@{},
   <-1.8mm,2.8mm>*{};<0.4mm,5.3mm>*{^1}**@{},
   <-2.8mm,2.9mm>*{};<-5.1mm,5.3mm>*{^3}**@{},
 \end{xy}}\Ea
 +
\Ba{c}\resizebox{8mm}{!}{\begin{xy}
 <0mm,0mm>*{\bu};<0mm,0mm>*{}**@{},
 <0mm,-0.49mm>*{};<0mm,-3.0mm>*{}**@{-},
 <0.49mm,0.49mm>*{};<1.9mm,1.9mm>*{}**@{-},
 <-0.5mm,0.5mm>*{};<-1.9mm,1.9mm>*{}**@{-},
 <-2.3mm,2.3mm>*{\bu};<-2.3mm,2.3mm>*{}**@{},
 <-1.8mm,2.8mm>*{};<0mm,4.9mm>*{}**@{-},
 <-2.8mm,2.9mm>*{};<-4.6mm,4.9mm>*{}**@{-},
   <0.49mm,0.49mm>*{};<2.7mm,2.3mm>*{^1}**@{},
   <-1.8mm,2.8mm>*{};<0.4mm,5.3mm>*{^3}**@{},
   <-2.8mm,2.9mm>*{};<-5.1mm,5.3mm>*{^2}**@{},
 \end{xy}}\Ea=0
 \ \ , \ \
\Ba{c}\resizebox{8.4mm}{!}{ \begin{xy}
 <0mm,0mm>*{\bu};<0mm,0mm>*{}**@{},
 <0mm,0.69mm>*{};<0mm,3.0mm>*{}**@{-},
 <0.39mm,-0.39mm>*{};<2.4mm,-2.4mm>*{}**@{-},
 <-0.35mm,-0.35mm>*{};<-1.9mm,-1.9mm>*{}**@{-},
 <-2.4mm,-2.4mm>*{\bu};<-2.4mm,-2.4mm>*{}**@{},
 <-2.0mm,-2.8mm>*{};<0mm,-4.9mm>*{}**@{-},
 <-2.8mm,-2.9mm>*{};<-4.7mm,-4.9mm>*{}**@{-},
    <0.39mm,-0.39mm>*{};<3.3mm,-4.0mm>*{^3}**@{},
    <-2.0mm,-2.8mm>*{};<0.5mm,-6.7mm>*{^2}**@{},
    <-2.8mm,-2.9mm>*{};<-5.2mm,-6.7mm>*{^1}**@{},
 \end{xy}}\Ea
 +
\Ba{c}\resizebox{8.4mm}{!}{ \begin{xy}
 <0mm,0mm>*{\bu};<0mm,0mm>*{}**@{},
 <0mm,0.69mm>*{};<0mm,3.0mm>*{}**@{-},
 <0.39mm,-0.39mm>*{};<2.4mm,-2.4mm>*{}**@{-},
 <-0.35mm,-0.35mm>*{};<-1.9mm,-1.9mm>*{}**@{-},
 <-2.4mm,-2.4mm>*{\bu};<-2.4mm,-2.4mm>*{}**@{},
 <-2.0mm,-2.8mm>*{};<0mm,-4.9mm>*{}**@{-},
 <-2.8mm,-2.9mm>*{};<-4.7mm,-4.9mm>*{}**@{-},
    <0.39mm,-0.39mm>*{};<3.3mm,-4.0mm>*{^2}**@{},
    <-2.0mm,-2.8mm>*{};<0.5mm,-6.7mm>*{^1}**@{},
    <-2.8mm,-2.9mm>*{};<-5.2mm,-6.7mm>*{^3}**@{},
 \end{xy}}\Ea
 +
\Ba{c}\resizebox{8.4mm}{!}{ \begin{xy}
 <0mm,0mm>*{\bu};<0mm,0mm>*{}**@{},
 <0mm,0.69mm>*{};<0mm,3.0mm>*{}**@{-},
 <0.39mm,-0.39mm>*{};<2.4mm,-2.4mm>*{}**@{-},
 <-0.35mm,-0.35mm>*{};<-1.9mm,-1.9mm>*{}**@{-},
 <-2.4mm,-2.4mm>*{\bu};<-2.4mm,-2.4mm>*{}**@{},
 <-2.0mm,-2.8mm>*{};<0mm,-4.9mm>*{}**@{-},
 <-2.8mm,-2.9mm>*{};<-4.7mm,-4.9mm>*{}**@{-},
    <0.39mm,-0.39mm>*{};<3.3mm,-4.0mm>*{^1}**@{},
    <-2.0mm,-2.8mm>*{};<0.5mm,-6.7mm>*{^3}**@{},
    <-2.8mm,-2.9mm>*{};<-5.2mm,-6.7mm>*{^2}**@{},
 \end{xy}}\Ea=0,
$$
$$
 \Ba{c}\resizebox{6mm}{!}{\begin{xy}
 <0mm,2.47mm>*{};<0mm,0.12mm>*{}**@{-},
 <0.5mm,3.5mm>*{};<2.2mm,5.2mm>*{}**@{-},
 <-0.48mm,3.48mm>*{};<-2.2mm,5.2mm>*{}**@{-},
 <0mm,3mm>*{\bu};<0mm,3mm>*{}**@{},
  <0mm,-0.8mm>*{\bu};<0mm,-0.8mm>*{}**@{},
<-0.39mm,-1.2mm>*{};<-2.2mm,-3.5mm>*{}**@{-},
 <0.39mm,-1.2mm>*{};<2.2mm,-3.5mm>*{}**@{-},
     <0.5mm,3.5mm>*{};<2.8mm,5.7mm>*{^2}**@{},
     <-0.48mm,3.48mm>*{};<-2.8mm,5.7mm>*{^1}**@{},
   <0mm,-0.8mm>*{};<-2.7mm,-5.2mm>*{^1}**@{},
   <0mm,-0.8mm>*{};<2.7mm,-5.2mm>*{^2}**@{},
\end{xy}}\Ea
+(-1)^{pq+p+q}\left(
\Ba{c}\resizebox{7mm}{!}{\begin{xy}
 <0mm,-1.3mm>*{};<0mm,-3.5mm>*{}**@{-},
 <0.38mm,-0.2mm>*{};<2.0mm,2.0mm>*{}**@{-},
 <-0.38mm,-0.2mm>*{};<-2.2mm,2.2mm>*{}**@{-},
<0mm,-0.8mm>*{\bu};<0mm,0.8mm>*{}**@{},
 <2.4mm,2.4mm>*{\bu};<2.4mm,2.4mm>*{}**@{},
 <2.77mm,2.0mm>*{};<4.4mm,-0.8mm>*{}**@{-},
 <2.4mm,3mm>*{};<2.4mm,5.2mm>*{}**@{-},
     <0mm,-1.3mm>*{};<0mm,-5.3mm>*{^1}**@{},
     <2.5mm,2.3mm>*{};<5.1mm,-2.6mm>*{^2}**@{},
    <2.4mm,2.5mm>*{};<2.4mm,5.7mm>*{^2}**@{},
    <-0.38mm,-0.2mm>*{};<-2.8mm,2.5mm>*{^1}**@{},
    \end{xy}}\Ea
  + (-1)^{q}
\Ba{c}\resizebox{7mm}{!}{\begin{xy}
 <0mm,-1.3mm>*{};<0mm,-3.5mm>*{}**@{-},
 <0.38mm,-0.2mm>*{};<2.0mm,2.0mm>*{}**@{-},
 <-0.38mm,-0.2mm>*{};<-2.2mm,2.2mm>*{}**@{-},
<0mm,-0.8mm>*{\bu};<0mm,0.8mm>*{}**@{},
 <2.4mm,2.4mm>*{\bu};<2.4mm,2.4mm>*{}**@{},
 <2.77mm,2.0mm>*{};<4.4mm,-0.8mm>*{}**@{-},
 <2.4mm,3mm>*{};<2.4mm,5.2mm>*{}**@{-},
     <0mm,-1.3mm>*{};<0mm,-5.3mm>*{^2}**@{},
     <2.5mm,2.3mm>*{};<5.1mm,-2.6mm>*{^1}**@{},
    <2.4mm,2.5mm>*{};<2.4mm,5.7mm>*{^2}**@{},
    <-0.38mm,-0.2mm>*{};<-2.8mm,2.5mm>*{^1}**@{},
    \end{xy}}\Ea
  + (-1)^{p+q}
\Ba{c}\resizebox{7mm}{!}{\begin{xy}
 <0mm,-1.3mm>*{};<0mm,-3.5mm>*{}**@{-},
 <0.38mm,-0.2mm>*{};<2.0mm,2.0mm>*{}**@{-},
 <-0.38mm,-0.2mm>*{};<-2.2mm,2.2mm>*{}**@{-},
<0mm,-0.8mm>*{\bu};<0mm,0.8mm>*{}**@{},
 <2.4mm,2.4mm>*{\bu};<2.4mm,2.4mm>*{}**@{},
 <2.77mm,2.0mm>*{};<4.4mm,-0.8mm>*{}**@{-},
 <2.4mm,3mm>*{};<2.4mm,5.2mm>*{}**@{-},
     <0mm,-1.3mm>*{};<0mm,-5.3mm>*{^2}**@{},
     <2.5mm,2.3mm>*{};<5.1mm,-2.6mm>*{^1}**@{},
    <2.4mm,2.5mm>*{};<2.4mm,5.7mm>*{^1}**@{},
    <-0.38mm,-0.2mm>*{};<-2.8mm,2.5mm>*{^2}**@{},
    \end{xy}}\Ea
 + (-1)^{p}
\Ba{c}\resizebox{7mm}{!}{\begin{xy}
 <0mm,-1.3mm>*{};<0mm,-3.5mm>*{}**@{-},
 <0.38mm,-0.2mm>*{};<2.0mm,2.0mm>*{}**@{-},
 <-0.38mm,-0.2mm>*{};<-2.2mm,2.2mm>*{}**@{-},
<0mm,-0.8mm>*{\bu};<0mm,0.8mm>*{}**@{},
 <2.4mm,2.4mm>*{\bu};<2.4mm,2.4mm>*{}**@{},
 <2.77mm,2.0mm>*{};<4.4mm,-0.8mm>*{}**@{-},
 <2.4mm,3mm>*{};<2.4mm,5.2mm>*{}**@{-},
     <0mm,-1.3mm>*{};<0mm,-5.3mm>*{^1}**@{},
     <2.5mm,2.3mm>*{};<5.1mm,-2.6mm>*{^2}**@{},
    <2.4mm,2.5mm>*{};<2.4mm,5.7mm>*{^1}**@{},
    <-0.38mm,-0.2mm>*{};<-2.8mm,2.5mm>*{^2}**@{},
    \end{xy}}\Ea\right)=0,
$$
which represent the Jacobi,  co-Jacobi and the Drinfeld identities.
To understand the deformation theory  of the properad $\LB_{p,q}$ one needs its minimal resolution $\HoLB_{p,q}$ which is a dg free properad generated by the  $(m,n)$-corollas
$$
\Ba{c}\resizebox{19mm}{!}{\begin{xy}
 <0mm,0mm>*{\bu};<0mm,0mm>*{}**@{},
 <-0.6mm,0.44mm>*{};<-8mm,5mm>*{}**@{-},
 <-0.4mm,0.7mm>*{};<-4.5mm,5mm>*{}**@{-},
 <0mm,0mm>*{};<1mm,5mm>*{\ldots}**@{},
 <0.4mm,0.7mm>*{};<4.5mm,5mm>*{}**@{-},
 <0.6mm,0.44mm>*{};<8mm,5mm>*{}**@{-},
   <0mm,0mm>*{};<-10.5mm,5.9mm>*{^{\sigma(1)}}**@{},
   <0mm,0mm>*{};<-4mm,5.9mm>*{^{\sigma(2)}}**@{},
   <0mm,0mm>*{};<10.0mm,5.9mm>*{^{\sigma(m)}}**@{},
 <-0.6mm,-0.44mm>*{};<-8mm,-5mm>*{}**@{-},
 <-0.4mm,-0.7mm>*{};<-4.5mm,-5mm>*{}**@{-},
 <0mm,0mm>*{};<1mm,-5mm>*{\ldots}**@{},
 <0.4mm,-0.7mm>*{};<4.5mm,-5mm>*{}**@{-},
 <0.6mm,-0.44mm>*{};<8mm,-5mm>*{}**@{-},
   <0mm,0mm>*{};<-10.5mm,-6.9mm>*{^{\tau(1)}}**@{},
   <0mm,0mm>*{};<-4mm,-6.9mm>*{^{\tau(2)}}**@{},
   <0mm,0mm>*{};<10.0mm,-6.9mm>*{^{\tau(n)}}**@{},
 \end{xy}}\Ea\hspace{-2mm}
=(-1)^{p|\sigma|+q|\tau|}
\Ba{c}\resizebox{14mm}{!}{\begin{xy}
 <0mm,0mm>*{\bu};<0mm,0mm>*{}**@{},
 <-0.6mm,0.44mm>*{};<-8mm,5mm>*{}**@{-},
 <-0.4mm,0.7mm>*{};<-4.5mm,5mm>*{}**@{-},
 <0mm,0mm>*{};<-1mm,5mm>*{\ldots}**@{},
 <0.4mm,0.7mm>*{};<4.5mm,5mm>*{}**@{-},
 <0.6mm,0.44mm>*{};<8mm,5mm>*{}**@{-},
   <0mm,0mm>*{};<-8.5mm,5.5mm>*{^1}**@{},
   <0mm,0mm>*{};<-5mm,5.5mm>*{^2}**@{},
   <0mm,0mm>*{};<4.5mm,5.5mm>*{^{m\hspace{-0.5mm}-\hspace{-0.5mm}1}}**@{},
   <0mm,0mm>*{};<9.0mm,5.5mm>*{^m}**@{},
 <-0.6mm,-0.44mm>*{};<-8mm,-5mm>*{}**@{-},
 <-0.4mm,-0.7mm>*{};<-4.5mm,-5mm>*{}**@{-},
 <0mm,0mm>*{};<-1mm,-5mm>*{\ldots}**@{},
 <0.4mm,-0.7mm>*{};<4.5mm,-5mm>*{}**@{-},
 <0.6mm,-0.44mm>*{};<8mm,-5mm>*{}**@{-},
   <0mm,0mm>*{};<-8.5mm,-6.9mm>*{^1}**@{},
   <0mm,0mm>*{};<-5mm,-6.9mm>*{^2}**@{},
   <0mm,0mm>*{};<4.5mm,-6.9mm>*{^{n\hspace{-0.5mm}-\hspace{-0.5mm}1}}**@{},
   <0mm,0mm>*{};<9.0mm,-6.9mm>*{^n}**@{},
 \end{xy}}\Ea
 \ \ \forall \sigma\in \bS_m \ \ \forall\tau\in \bS_n, \ m,n\geq 1, m+n\geq 3,
$$
of degrees $1+p(1-m) +q(1-n)$ on which the differential acts by splitting them as follows
$$
\delta
\Ba{c}\resizebox{14mm}{!}{\begin{xy}
 <0mm,0mm>*{\bu};<0mm,0mm>*{}**@{},
 <-0.6mm,0.44mm>*{};<-8mm,5mm>*{}**@{-},
 <-0.4mm,0.7mm>*{};<-4.5mm,5mm>*{}**@{-},
 <0mm,0mm>*{};<-1mm,5mm>*{\ldots}**@{},
 <0.4mm,0.7mm>*{};<4.5mm,5mm>*{}**@{-},
 <0.6mm,0.44mm>*{};<8mm,5mm>*{}**@{-},
   <0mm,0mm>*{};<-8.5mm,5.5mm>*{^1}**@{},
   <0mm,0mm>*{};<-5mm,5.5mm>*{^2}**@{},
   <0mm,0mm>*{};<4.5mm,5.5mm>*{^{m\hspace{-0.5mm}-\hspace{-0.5mm}1}}**@{},
   <0mm,0mm>*{};<9.0mm,5.5mm>*{^m}**@{},
 <-0.6mm,-0.44
 mm>*{};<-8mm,-5mm>*{}**@{-},
 <-0.4mm,-0.7mm>*{};<-4.5mm,-5mm>*{}**@{-},
 <0mm,0mm>*{};<-1mm,-5mm>*{\ldots}**@{},
 <0.4mm,-0.7mm>*{};<4.5mm,-5mm>*{}**@{-},
 <0.6mm,-0.44mm>*{};<8mm,-5mm>*{}**@{-},
   <0mm,0mm>*{};<-8.5mm,-6.9mm>*{^1}**@{},
   <0mm,0mm>*{};<-5mm,-6.9mm>*{^2}**@{},
   <0mm,0mm>*{};<4.5mm,-6.9mm>*{^{n\hspace{-0.5mm}-\hspace{-0.5mm}1}}**@{},
   <0mm,0mm>*{};<9.0mm,-6.9mm>*{^n}**@{},
 \end{xy}}\Ea
\ \ = \ \
 \sum_{[m]=I_1\sqcup I_2\atop
 [n]=J_1\sqcup J_2}
\pm
\Ba{c}\resizebox{22mm}{!}{ \begin{xy}
 <0mm,0mm>*{\bu};<0mm,0mm>*{}**@{},
 <-0.6mm,0.44mm>*{};<-8mm,5mm>*{}**@{-},
 <-0.4mm,0.7mm>*{};<-4.5mm,5mm>*{}**@{-},
 <0mm,0mm>*{};<0mm,5mm>*{\ldots}**@{},
 <0.4mm,0.7mm>*{};<4.5mm,5mm>*{}**@{-},
 <0.6mm,0.44mm>*{};<12.4mm,4.8mm>*{}**@{-},
     <0mm,0mm>*{};<-2mm,7mm>*{\overbrace{\ \ \ \ \ \ \ \ \ \ \ \ }}**@{},
     <0mm,0mm>*{};<-2mm,9mm>*{^{I_1}}**@{},
 <-0.6mm,-0.44mm>*{};<-8mm,-5mm>*{}**@{-},
 <-0.4mm,-0.7mm>*{};<-4.5mm,-5mm>*{}**@{-},
 <0mm,0mm>*{};<-1mm,-5mm>*{\ldots}**@{},
 <0.4mm,-0.7mm>*{};<4.5mm,-5mm>*{}**@{-},
 <0.6mm,-0.44mm>*{};<8mm,-5mm>*{}**@{-},
      <0mm,0mm>*{};<0mm,-7mm>*{\underbrace{\ \ \ \ \ \ \ \ \ \ \ \ \ \ \
      }}**@{},
      <0mm,0mm>*{};<0mm,-10.6mm>*{_{J_1}}**@{},
 <13mm,5mm>*{};<13mm,5mm>*{\bu}**@{},
 <12.6mm,5.44mm>*{};<5mm,10mm>*{}**@{-},
 <12.6mm,5.7mm>*{};<8.5mm,10mm>*{}**@{-},
 <13mm,5mm>*{};<13mm,10mm>*{\ldots}**@{},
 <13.4mm,5.7mm>*{};<16.5mm,10mm>*{}**@{-},
 <13.6mm,5.44mm>*{};<20mm,10mm>*{}**@{-},
      <13mm,5mm>*{};<13mm,12mm>*{\overbrace{\ \ \ \ \ \ \ \ \ \ \ \ \ \ }}**@{},
      <13mm,5mm>*{};<13mm,14mm>*{^{I_2}}**@{},
 <12.4mm,4.3mm>*{};<8mm,0mm>*{}**@{-},
 <12.6mm,4.3mm>*{};<12mm,0mm>*{\ldots}**@{},
 <13.4mm,4.5mm>*{};<16.5mm,0mm>*{}**@{-},
 <13.6mm,4.8mm>*{};<20mm,0mm>*{}**@{-},
     <13mm,5mm>*{};<14.3mm,-2mm>*{\underbrace{\ \ \ \ \ \ \ \ \ \ \ }}**@{},
     <13mm,5mm>*{};<14.3mm,-4.5mm>*{_{J_2}}**@{},
 \end{xy}}\Ea
$$
The natural projection
\Beq\label{2: p from HoLB toLB}
\pi: \HoLB_{p,q} \lon \LB_{p,q}
\Eeq
which sends to zero all non-trivalent generators, is a  quasi-isomorphism.

\sip

The (genus completed)
derivation complex
of the properad  $\HoLB_{p,q}$ has been studied in detail in \cite{MW1} where an explicit quasi-isomorphism of complexes
\Beq\label{2: F from OGC to Der(HoLB)}
F: \GC_{p+q+1}^{\mathsf{or}, \geq 2} \oplus \K \lon \Der(\HoLB_{p,q})
\Eeq
has been constructed by sending an oriented graph $\Ga$ to a derivation $F(\Ga)$ whose value on the
$(m,n)$-generator of $\HoLB_{p,q}$ is given explicitly by
given explicitly by
$$
F(\Ga): \Ba{c}\resizebox{12mm}{!}{\begin{xy}
 <0mm,0mm>*{\circ};<0mm,0mm>*{}**@{},
 <-0.6mm,0.44mm>*{};<-8mm,5mm>*{}**@{-},
 <-0.4mm,0.7mm>*{};<-4.5mm,5mm>*{}**@{-},
 <0mm,0mm>*{};<-1mm,5mm>*{\ldots}**@{},
 <0.4mm,0.7mm>*{};<4.5mm,5mm>*{}**@{-},
 <0.6mm,0.44mm>*{};<8mm,5mm>*{}**@{-},
   <0mm,0mm>*{};<-8.5mm,5.5mm>*{^1}**@{},
   <0mm,0mm>*{};<-5mm,5.5mm>*{^2}**@{},
   <0mm,0mm>*{};<9.0mm,5.5mm>*{^m}**@{},
 <-0.6mm,-0.44mm>*{};<-8mm,-5mm>*{}**@{-},
 <-0.4mm,-0.7mm>*{};<-4.5mm,-5mm>*{}**@{-},
 <0mm,0mm>*{};<-1mm,-5mm>*{\ldots}**@{},
 <0.4mm,-0.7mm>*{};<4.5mm,-5mm>*{}**@{-},
 <0.6mm,-0.44mm>*{};<8mm,-5mm>*{}**@{-},
   <0mm,0mm>*{};<-8.5mm,-6.9mm>*{^1}**@{},
   <0mm,0mm>*{};<-5mm,-6.9mm>*{^2}**@{},
   <0mm,0mm>*{};<9.0mm,-6.9mm>*{^n}**@{},
 \end{xy}}\Ea
 \lon
 \sum_{s:[n]\rar V(\Ga)\atop \hat{s}:[m]\rar V(\Ga)}  \Ba{c}\resizebox{9mm}{!}  {\xy
 (-6,7)*{^1},
(-3,7)*{^2},
(2.5,7)*{},
(7,7)*{^m},
(-3,-8)*{_2},
(3,-6)*{},
(7,-8)*{_n},
(-6,-8)*{_1},
(0,4.5)*+{...},
(0,-4.5)*+{...},
(0,0)*+{\Ga}="o",
(-6,6)*{}="1",
(-3,6)*{}="2",
(3,6)*{}="3",
(6,6)*{}="4",
(-3,-6)*{}="5",
(3,-6)*{}="6",
(6,-6)*{}="7",
(-6,-6)*{}="8",
\ar @{-} "o";"1" <0pt>
\ar @{-} "o";"2" <0pt>
\ar @{-} "o";"3" <0pt>
\ar @{-} "o";"4" <0pt>
\ar @{-} "o";"5" <0pt>
\ar @{-} "o";"6" <0pt>
\ar @{-} "o";"7" <0pt>
\ar @{-} "o";"8" <0pt>
\endxy}\Ea
$$
 where the sum is taken over all ways of attaching the incoming and outgoing legs to the graph $\Ga$, and one sets to zero every resulting graph if it contains a vertex with valency $<3$ or
   with no at least one incoming  or at least one outgoing edge.

   \sip

   As there is  an isomorphism of Lie algebras \cite{Wi1,Wi2}
$$
H^0(\GC_3^{\mathsf{or, \geq 2}})\oplus \K \simeq H^0(\GC_2^{\mathsf{\geq 2}}) \oplus \K=\grt_1\oplus \K \simeq \grt
$$
where the direct sum symbol is applied in the category of graded vector spaces (not in the category of Lie algebras), one concludes \cite{MW1} that the mysterious Grothendieck-Teichm\"uller group $GRT$ introduced by V.\ Drinfeld in \cite{Dr} acts faithfully on the genus completion of the properad $\LB_{1,1}$ of ordinary Lie bialgebras as homotopy non-trivial automorphisms.

\subsection{Cohomology of oriented graph complexes in terms of bivalent and trivalent graphs} Using the epimorphism (\ref{2: p from HoLB toLB}) one can consider the derivation complex $\Der(\HoLB_{p,q}\rar \LB_{p,q})$  of the properad $\HoLB_{p,q}$ with values in $\LB_{p,q}$; up to degree shift by $1$ this complex is identical to the deformation
complex $\Def(\HoLB_{p,q}\rar \LB_{p,q})$
 of the epimorphism $\pi$ (see \cite{MV,MW1} for more details). As the map  (\ref{2: p from HoLB toLB}) is a quasi-isomorphism,
the induced morphism of complexes
$$
P: \Der(\HoLB_{p,q}) \lon \Der(\HoLB_{p,q}\rar \LB_{p,q})
$$
is also a quasi-isomorphism. The latter derivation complex is generated by {\em trivalent graphs with in- and out-legs (``hairs")}\, as the generating corollas of $\LB_{p,q}$ are of this type.
Composing the latter with the morphism $F$ in (\ref{2: F from OGC to Der(HoLB)}) one concludes that graphs in $\GC^{\mathsf{or, \geq 2}}_{p+q+1}$ with valencies of vertices $\geq 4$ do not contribute into the cohomology groups. More precisely,
let ${\GC}^{\mathsf{3,\geq 4}}_{d+1}$ be the differential closure of the linear
subspace of
$\GC^{\mathsf{or,\geq 2}}_{d+1}$ spanned by graphs with at least one vertex of valency
$\geq 4$ or with at least one trivalent source or target, and let  $\GC_{d+1}^{\mathsf{lieb}}$ be
quotient complex,
 $$
0\lon \GC^{\mathsf{3,\geq 4}}_{d+1} \lon \GC^{\mathsf{or, \geq 2}}_{d+1} \stackrel{\pi^{\mathsf{or}}}{\lon} \GC_{d+1}^{\mathsf{lieb}}\lon 0
$$
 which is generated by oriented graphs with vertices of types $(2,1)$, $(1,2)$,
 $(2,0)$ and $(0,2)$ modulo $\LB$-type relations. The induced differential in $\GC_{d+1}^{\mathsf{lieb}}$  acts only on $(1,2)$ and $(2,1)$
vertices creating new bivalent targets and sources.
There is a commutative diagram of short exact sequences of complexes,
$$
\Ba{c}\resizebox{120mm}{!}{
\xymatrix{
0\ar[r]& \mbox{$\GC
_{p+q+1}^{3,\geq 4}$}\ar[d]\ar[r] &
\mbox{$\GC^{\mathsf{or},\geq 2}_{p+q+1}$} \ar[r]\ar[d] &
\mbox{$\GC_{p+q+1}^{\mathsf{lieb}}$}\ar[d]\ar[r]  & 0\\
0\ar[r]& \ker P\ar[r]      &
\Der(\HoLB_{p,q}) \ar[r]    &  \Der(\HoLB_{p,q}\rar \LB_{p,q}) \ar[r]       & 0\\
}}\Ea
$$
As the complex $\ker P$ is acyclic, this diagram implies that the induced morphism of cohomology groups
$$
 H^\bu(\GC^{\mathsf{or},\geq 2}_{d+1})\lon H^\bu(\GC_{d+1}^{\mathsf{lieb}})
$$
is an {\em injection}\, for any $d\in \Z$. It is easy to improve this observation to a stronger statement as follows.

 \subsubsection{\bf Proposition} {\em The canonical projection
$
\pi^{\mathsf{or}}: \GC^{\mathsf{or}, \geq 2}_{d+1} \rar {\GC}_{d+1}^{\mathsf{lieb}}
$
is a quasi-isomorphism for any $d\in \Z$.}

\begin{proof}
Consider a filtration of both sides of the epimorphism $\pi^{\mathsf{or}}$ by the total number of sources and targets. One gets an induced morphism of the associated graded complexes
$$
gr(\pi^{\mathsf{or}}): (gr\GC^{\mathsf{or}, \geq 2}_{d+1}, \delta_0) \lon (gr{\GC}_{d+1}^{\mathsf{lieb}},0)
$$
As the number targets and sources of graphs in both complexes just above are preserved by the induced differentials,  one can consider
a version of the above diagram
$$
\hat{gr}(\pi^{\mathsf{or}}): (\wh{gr}\GC^{\mathsf{or}, \geq 2}_{d+1}, \delta_0) \lon
(\wh{gr}{\GC}_{d+1}^{\mathsf{lieb}},0)
$$
in which sources and targets of the generating graphs are distinguished, say, by attaching to each target (resp. source) a labelled out-going (resp., ingoing) hair,
$$
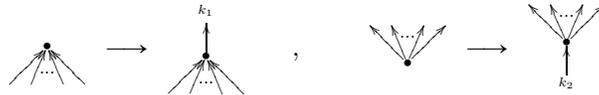
\Ba{c}
\resizebox{10mm}{!}{\xy
(0,-4.9)*+{...},
(0,0)*{\bu}="o",
(0,6)*{}="1",
(-3,-7)*{}="5",
(3,-7)*{}="6",
(7,-7)*{}="7",
(-7,-7)*{}="8",
%
%
\ar @{<-} "o";"5" <0pt>
\ar @{<-} "o";"6" <0pt>
\ar @{<-} "o";"7" <0pt>
\ar @{<-} "o";"8" <0pt>
\endxy}\Ea
\lon
\Ba{c}
\resizebox{10mm}{!}{\xy
(0,8)*{^{k_1}},
(0,-4.9)*+{...},
(0,0)*{\bu}="o",
(0,6)*{}="1",
(-3,-7)*{}="5",
(3,-7)*{}="6",
(7,-7)*{}="7",
(-7,-7)*{}="8",
\ar @{->} "o";"1" <0pt>
\ar @{<-} "o";"5" <0pt>
\ar @{<-} "o";"6" <0pt>
\ar @{<-} "o";"7" <0pt>
\ar @{<-} "o";"8" <0pt>
\endxy}\Ea
\ \ \ \ , \ \ \ \ \ \
\Ba{c}
\resizebox{10mm}{!}{\xy
(0,4.9)*+{...},
(0,0)*{\bu}="o",
(0,6)*{}="1",
(-3,7)*{}="5",
(3,7)*{}="6",
(7,7)*{}="7",
(-7,7)*{}="8",
%
%
\ar @{->} "o";"5" <0pt>
\ar @{->} "o";"6" <0pt>
\ar @{->} "o";"7" <0pt>
\ar @{->} "o";"8" <0pt>
\endxy}\Ea
\lon
\Ba{c}
\resizebox{10mm}{!}{\xy
(0,-8)*{^{k_2}},
(0,4.9)*+{...},
(0,0)*{\bu}="o",
(0,-6)*{}="1",
(-3,7)*{}="5",
(3,7)*{}="6",
(7,7)*{}="7",
(-7,7)*{}="8",
\ar @{<-} "o";"1" <0pt>
\ar @{->} "o";"5" <0pt>
\ar @{->} "o";"6" <0pt>
\ar @{->} "o";"7" <0pt>
\ar @{->} "o";"8" <0pt>
\endxy}\Ea
$$
Now one can identify directly the complex $(\wh{gr}\GC^{\mathsf{or, \geq 2}}_{d+1}, \delta_0)$ with the complex $\HoLB_{0,d}$ while the complex $(\wh{gr}{\GC}_d^{lieb},0)$ with its cohomology
$\LB_{0,d}$. As the projection (\ref{2: p from HoLB toLB}) is a quasi-isomorphism,  the map  $\wh{gr}(\pi^\mathsf{or})$  is a quasi-isomorphism as well implying,
by the Maschke theorem,  that the map  $gr(\pi^{\mathsf{or}})$ is a quasi-isomorphism which in turn implies that the map $\pi^\mathsf{or}$ is a quasi-isomorphism.
\end{proof}

There is another ``small" description of the cohomology of $H^\bu(\OGCd)\simeq H^\bu(\GCd)$
in terms of the cohomology of a complex generated by  {\em solely trivalent graphs}\, which we describe in the next subsection.

 \subsection{A purely trivalent incarnation of the Kontsevich graph complex}
Let $\overline{\GC}^{\mathsf{or,\geq 4}}_{d+1}$ be the differential closure of the subspace in $\OOGCd$
generated by graphs having at least one vertex of valency $\geq 4$, and consider the quotient complex  $\oGCD$,
$$
0\lon \overline{\GC}^{or,\geq 4}_{d+1} \lon \OOGCd \stackrel{s}{\lon} \oGCD \lon 0.
$$
By its very construction, $\oGCD$ is a graded vector space
generated  by {\em trivalent} graphs $\Ga$ with two types of edges,
solid edges
\ { $\xy
 (0,0)*{\bullet}="a",
(7,0)*{\bu}="b",
\ar @{->} "a";"b" <0pt>
\endxy$} equipped with a fixed direction, and dotted ones whose direction can be flipped
{$
 
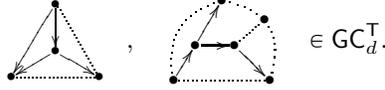
\xy
 (0,0)*{\bullet}="a",
(7,0)*{\bu}="b",
\ar @{.>} "a";"b" <0pt>
\endxy = (-1)^{d} \xy
 (0,0)*{\bullet}="a",
(7,0)*{\bu}="b",
\ar @{<.} "a";"b" <0pt>
\endxy$}, e.g.
 $$
\Ba{c}\resizebox{14mm}{!}{
\xy
 (0,0)*{\bullet}="a",
(0,8)*{\bullet}="b",
(-7.5,-4.5)*{\bu}="c",
(7.5,-4.5)*{\bu}="d",
\ar @{<-} "a";"b" <0pt>
\ar @{->} "a";"c" <0pt>
\ar @{->} "b";"c" <0pt>
\ar @{.} "d";"c" <0pt>
\ar @{.} "b";"d" <0pt>
\ar @{<-} "d";"a" <0pt>
\endxy}
\Ea
 \ , \
  \Ba{c}\resizebox{16mm}{!}{
\xy
 (0,8)*{\bullet}="a",
(-7.5,-4.5)*{\bullet}="c",
(7.5,-4.5)*{\bu}="d",
(6.7,5)*{\bu}="r",
(2,1)*{\bu}="o",
(-4.1,1)*{\bu}="l",
\ar@/^-0.8pc/@{.}"a";"c" <0pt>
\ar@/^0.99pc/@{.}"a";"d" <0pt>
\ar @{<-} "l";"c" <0pt>
\ar @{->} "l";"a" <0pt>
\ar @{->} "l";"o" <0pt>
\ar @{.} "d";"c" <0pt>
\ar @{->} "o";"d" <0pt>
\ar @{.} "o";"r" <0pt>
\endxy}
\Ea \ \in \oGCD.
$$
Solid edges of $\Ga$ are not allowed to form a closed path of directed edges (a so called {\em solid wheel}), and the
cohomological degree is given by the formula
 $$
 |\Ga|=(d+1)(\# V(\Ga)-1) +(1-d)\# E_{dot}(\Ga) - d\#E_{sol}(\Ga).
 $$
Such trivalent graphs are required to satisfy an IHX-type relation for any its {\em solid}\, edge whose contraction does not create solid wheels,
$$
\left(\hspace{-2mm}
\Ba{c}\resizebox{12mm}{!}{
\xy
 (-4,0)*{\bullet}="a",
(+4,0)*{\bu}="b",
 (-6,11)*+{^1}="1",
(-6,-10)*+{_2}="2",
 (6,11)*+{^3}="3",
(6,-10)*+{_4}="4",
\ar @{->} "a";"b" <0pt>
\ar @{--} "a";"1" <0pt>
\ar @{--} "a";"2" <0pt>
\ar @{--} "b";"3" <0pt>
\ar @{--} "b";"4" <0pt>
\endxy}\Ea
\hspace{-2mm}
+
\hspace{-2mm}
\Ba{c}\resizebox{12mm}{!}{
\xy
 (-4,0)*{\bullet}="a",
(+4,0)*{\bu}="b",
 (-6,11)*+{^1}="1",
(-6,-10)*+{_3}="2",
 (6,11)*+{^2}="3",
(6,-10)*+{_4}="4",
\ar @{->} "a";"b" <0pt>
\ar @{--} "a";"1" <0pt>
\ar @{--} "a";"2" <0pt>
\ar @{--} "b";"3" <0pt>
\ar @{--} "b";"4" <0pt>
\endxy}\Ea
\hspace{-2mm}
+
\hspace{-2mm}
\Ba{c}\resizebox{11mm}{!}{
\xy
 (-4,0)*{\bullet}="a",
(+4,0)*{\bu}="b",
 (-6,11)*+{^1}="1",
(-6,-10)*+{_4}="2",
 (6,11)*+{^3}="3",
(6,-10)*+{_2}="4",
\ar @{->} "a";"b" <0pt>
\ar @{--} "a";"1" <0pt>
\ar @{--} "a";"2" <0pt>
\ar @{--} "b";"3" <0pt>
\ar @{--} "b";"4" <0pt>
\endxy}\Ea
\hspace{-2mm}
\right)
-(-1)^d
\left(\hspace{-2mm}
\Ba{c}\resizebox{12mm}{!}{
\xy
 (-4,0)*{\bullet}="a",
(+4,0)*{\bu}="b",
 (-6,11)*+{^1}="1",
(-6,-10)*+{_2}="2",
 (6,11)*+{^3}="3",
(6,-10)*+{_4}="4",
\ar @{<-} "a";"b" <0pt>
\ar @{--} "a";"1" <0pt>
\ar @{--} "a";"2" <0pt>
\ar @{--} "b";"3" <0pt>
\ar @{--} "b";"4" <0pt>
\endxy}\Ea
\hspace{-2mm}
+
\hspace{-2mm}
\Ba{c}\resizebox{12mm}{!}{
\xy
 (-4,0)*{\bullet}="a",
(+4,0)*{\bu}="b",
 (-6,11)*+{^1}="1",
(-6,-10)*+{_3}="2",
 (6,11)*+{^2}="3",
(6,-10)*+{_4}="4",
\ar @{<-} "a";"b" <0pt>
\ar @{--} "a";"1" <0pt>
\ar @{--} "a";"2" <0pt>
\ar @{--} "b";"3" <0pt>
\ar @{--} "b";"4" <0pt>
\endxy}\Ea
\hspace{-2mm}
+
\hspace{-2mm}
\Ba{c}\resizebox{11mm}{!}{
\xy
 (-4,0)*{\bullet}="a",
(+4,0)*{\bu}="b",
 (-6,11)*+{^1}="1",
(-6,-10)*+{_4}="2",
 (6,11)*+{^3}="3",
(6,-10)*+{_2}="4",
\ar @{<-} "a";"b" <0pt>
\ar @{--} "a";"1" <0pt>
\ar @{--} "a";"2" <0pt>
\ar @{--} "b";"3" <0pt>
\ar @{--} "b";"4" <0pt>
\endxy}\Ea
\hspace{-2mm}
\right)
=0,
$$
where every dashed edge stands for an arbitrary solid (of any direction) or dotted edge. This relation originates from the part $\delta_\bu$ of the full differential in $\OOGCd$ applied to a 4-valent vertex.
\sip

The induced  differential on $\oGCD$ acts only on solid edges of the generators $\Ga$ by making them dotted, i.e.
$$
\delta\Ga:= \sum_{e\in E_{sol}(\Ga)} \delta_e \Ga \ \ \text{where}\ \
\delta_e:   \Ba{c}\resizebox{11mm}{!}{
\xy
(0,2)*{^{v_1}},
(8,2)*{^{v_2}},
 (0,0)*{\bullet}="a",
(8,0)*{\bu}="b",
\ar @{->} "a";"b" <0pt>
\endxy}\Ea
\rightsquigarrow
 \Ba{c}\resizebox{11mm}{!}{
\xy
(0,2)*{^{v_1}},
(8,2)*{^{v_2}},
 (0,0)*{\bullet}="a",
(8,0)*{\bu}="b",
\ar @{.>} "a";"b" <0pt>
\endxy}\Ea.
$$
Note that $\delta$ respects the IHX-relation as making the solid edge in its center into a dotted one makes the whole combination of graphs vanishing identically.

\subsubsection{\bf Proposition}\label{3: Prop on hatOGC^3} {\em
The canonical projection
$s: \OOGCd {\lon} \oGCD$
is a quasi-isomorphism.
}



\begin{proof}
Consider a filtration of both sides of the morphism $s$
 by the number of dotted edges. The induced differential in the associated graded complex $gr \oGCD$ is trivial,
while the induced differential in $gr \OOGCd$ acts by splitting vertices
and creating a new solid edge. Since the set of dotted edges is preserved by the induced differential, we can consider a version of  $gr \OOGCd$
in which the dotted edges are distinguished and directed, i.e. the half-edges forming a dotted edge
are totally ordered. One can then cut every directed dotted edge into a pair of directed legs (or hairs), one is incoming and
one is outgoing, both being equipped with the same label as the one assigned to the original dotted edge.
The resulting complex can then be identified with direct summand $\oplus_{m\geq 1} p\HoLB_{0,d}(m,m)$ of  the degree shifted
 dg properad $p\HoLB_{0,d}=\{p\HoLB_{0,d}(m,n)\} $, the
minimal resolution of properad the so called {\em pseudo}-Lie bialgebras $p\LB_{0,d}$ (which is a version of properad Lie bialgebras $\LB_{0,d}$ in which trivalent
sources and targets are allowed, see \cite{G} for the detailed definition and the study of this properad).  We conclude that the cohomology of the latter complex is generated by trivalent graphs only, i.e.\ it coincides with the corresponding haired version  of $\oGCD$.
Hence the projection $s$ induces an isomorphism of the second pages of
the spectral sequence by the number of vertices. By the Comparison Theorem, the map $s$
is a quasi-isomorphism.
\end{proof}

\sip

Hence one has an isomorphism of cohomology groups
$$
H^\bu(\oGCD)\simeq H^\bu(\GCd),
$$
i.e\ the complex $\oGCD$ gives us an incarnation of the M.\ Kontsevich graph complex in terms of purely trivalent graphs.
It is worth noting that a graph $\Ga\in \oGCD$ with a fixed cohomological degree $|\Ga|$ and a fixed loop number $g$ must satisfy the following conditions,
$$
\# V(\Ga)=2g-2, \ \ \ \# E_{dot}(\Ga)=|\Ga|+3 +g(d-2) , \ \ \ \# E_{sol}(\Ga)=-|\Ga| + g(5-d) -6.
$$

\sip

The above trivalent representation of $\OOGCd$ suggests that the cohomology
$H^\bu(\GCd)$ of the original M.\ Kontsevich graph complex might be concentrated  in a quotient space generated by graphs with valencies of vertices bounded from above. There is a conjecture (stemming from the deformation quantization theory of quadratic Poisson structures) which  says that graphs with trivalent and four-valent vertices only might do the job.

\subsubsection{\bf Conjecture} {\em Let $\GCd^{\geq 5}$ be the
differential closure of the
subspace in $\GCd$ generated by graphs with at least one vertex of valency $\geq 5$. Then
the epimorphism $\pi$ in the short exact sequence
$$
0\lon \GCd^{\geq 5} \lon \GCd \stackrel{\pi}{\lon} \GCd^{\leq 4}\lon 0
$$
is a quasi-isomorphism}.

\sip

Simon Brun and Thomas Willwacher used numerical methods in \cite{BW} to show that
the projection $\GC_{g,d} \stackrel{\pi}{\rar} \GC_{g,d}^{\leq 4}$ is indeed a quasi-isomorphism for all $g\leq 10$.


\mip

{\large
\section{\bf Directed graph complex: its reduced versions}
}

\subsection{A reduced version of $\dGCD$  with two types of edges}\label{4: subsec on reduced dGC with 2 types of edges} Let $\rrdGCD$ stand for an extension of the complex $\rOGCd$ (see \S 2.5) in which the condition that solid edges never form closed paths of directed edges is dropped. Put another way,  $\rrdGCD$ is generated by graphs with all vertices {\em at least trivalent}\, and with two types of edges, solid edges
$
 \xy
 (0,0)*{\bullet}="a",
(7,0)*{\bu}="b",
\ar @{->} "a";"b" <0pt>
\endxy$  which have a fixed direction and degree $-d$   and dotted edges
 $\xy
 (0,0)*{\bullet}="a",
(7,0)*{\bu}="b",
\ar @{.>} "a";"b" <0pt>
\endxy.
$ which have degree $1-d$ and whose direction can be flipped as in (\ref{2: symmetry of dotted edges}). The differential in $\rrdGCD$ is given by formula (\ref{2: delta on hatOGCd}).

\subsubsection{\bf Proposition}\label{4: Prop on 2-edge version of dGC} {\em There is an explicit quasi-isomorphism of complexes
$$
\Ba{rccc}
s:& \rrdGCD &\lon & \dGCD    \\
 & \Ga & \lon & s(\Ga)
\Ea
$$
where $s(\Ga)$ is obtained from $\Ga$ by changing its edges according to the rule (cf.\ (\ref{2: map j from dooted edge to solid comb}))
$$
s: \left\{
\Ba{l}
 \Ba{c}\resizebox{11mm}{!}{
\xy
(0,2.5)*{^{v_1}},
(8,2.5)*{^{v_2}},
 (0,1)*{\bullet}="a",
(8,1)*{\bu}="b",
\ar @{->} "a";"b" <0pt>
\endxy}\Ea
  \rightsquigarrow
   \Ba{c}\resizebox{11mm}{!}{
\xy
(0,2.5)*{^{v_1}},
(8,2.5)*{^{v_2}},
 (0,1)*{\bullet}="a",
(8,1)*{\bu}="b",
\ar @{->} "a";"b" <0pt>
\endxy}\Ea
\vspace{2mm}
\\
\Ba{c}\resizebox{11mm}{!}{
\xy
(0,2.5)*{^{v_1}},
(8,2.5)*{^{v_2}},
 (0,1)*{\bullet}="a",
(8,1)*{\bu}="b",
\ar @{.>} "a";"b" <0pt>
\endxy}\Ea
 \rightsquigarrow
 \frac{1}{2}\left(
\Ba{c}\resizebox{15mm}{!}{  \xy
(0,2)*{^{v_1}},
(7,5)*{^{v}},
(14,2)*{^{v_2}},
 (0,0)*{\bullet}="a",
 (7,3)*{\bu}="c",
(14,0)*{\bu}="b",
\ar @{->} "a";"c" <0pt>
\ar @{->} "b";"c" <0pt>
\endxy}\Ea
-
\Ba{c}\resizebox{15mm}{!}{
\xy
(0,2)*{^{v_1}},
(7,-1)*{^{v}},
(14,2)*{^{v_2}},
 (0,0)*{\bullet}="a",
 (7,-3)*{\bu}="c",
(14,0)*{\bu}="b",
\ar @{<-} "a";"c" <0pt>
\ar @{<-} "b";"c" <0pt>
\endxy}\Ea
\right)
\Ea
\right.
$$
}

\begin{proof} Consider a filtration of both sides of the map $s$ by the number of vertices. The induced differential in the associated graded complex of the left (resp., right) hand side is $\delta'$ (resp., trivial). The cohomology of the former complex is generated by graphs with all directed edges (skew)symmetrized,   $\Ba{c}\resizebox{11mm}{!}{
\xy
(0,2.5)*{^{v_1}},
(8,2.5)*{^{v_2}},
 (0,1)*{\bullet}="a",
(8,1)*{\bu}="b",
\ar @{->} "a";"b" <0pt>
\endxy}\Ea -  (-1)^d \Ba{c}\resizebox{11mm}{!}{
\xy
(0,2.5)*{^{v_1}},
(8,2.5)*{^{v_2}},
 (0,1)*{\bullet}="a",
(8,1)*{\bu}="b",
\ar @{<-} "a";"b" <0pt>
\endxy}\Ea$. Hence the next pages of the spectral sequences give us an inclusion of complexes
$$
\Img (f) \lon \dGCD
$$
which is a quasi-isomorphism \cite{Wi1}; here $f$ is the map (\ref{1: GC to dGC}).
By the Comparison Theorem the map $s$ is a quasi-isomorphism.
\end{proof}

Let $\rTGCd$ be a subcomplex of $\rrdGCD$ spanned by graphs having at least one vertex (called a {\em target}) with no outgoing solid edges. Similarly one defines a subcomplex $\overline{\GC}_{d+1}^{\mathsf{s}}\subset \rrdGCD$ spanned by graphs having at least one {\em source}, that is  a vertex
with no ingoing solid edges. These sub-complexes give us reduced models for targeted and sourced graph complexes $\TGCd$ and $\SGCd$ respectively (see \S 2.4 for their definitions).

\subsubsection{\bf Theorem}\label{4: Theorem rTGCd versis TGCd} {\em For any $d\in \Z$ one has}
$$
H^\bu(\rTGCd)=H^\bu(\TGCd), \ \ \ H^\bu(\overline{\GC}_{d+1}^{\mathsf{s}})=H^\bu(\SGCd)
$$

This result may seem a bit surprising as the dotted edge in $\rrdGCD$ stands for a linear combination of a bivalent source {\em and}\, a bivalent target. We show its proof below after reminding in the next subsection another useful  reduced model for $\dGCD$ in which sources and targets are clearly distinguished.

\subsection{Marko \v Zivkovi\'c' reduced version of $\dGCD$ with four types of edges}
\label{4: subsec on MZ reduced complex rdGCD}
 Let $\rdGCD$ be a graded vector space generated by graphs $\Ga$  which have all vertices
 $v\in V(\Ga)$ at least trivalent and which can have edges of the following four types,
\Bi
\item[(i)] solid edges  $\Ba{c}\resizebox{11mm}{!}{
\xy
(0,2.5)*{^{v_1}},
(8,2.5)*{^{v_2}},
 (0,1)*{\bullet}="a",
(8,1)*{\bu}="b",
\ar @{->} "a";"b" <0pt>
\endxy}\Ea$ equipped with a fixed direction;
\item[(ii)] two types of dotted edges, dotted $s$-edges  $\Ba{c}\resizebox{11mm}{!}{
\xy
(0,2.5)*{^{v_1}},
(8,2.5)*{^{v_2}},
(4,2)*{^{s}},
 (0,1)*{\bullet}="a",
(8,1)*{\bu}="b",
\ar @{.>} "a";"b" <0pt>
\endxy}\Ea$   and dotted $t$-edges  $\Ba{c}\resizebox{11mm}{!}{
\xy
(0,2.5)*{^{v_1}},
(8,2.5)*{^{v_2}},
(4,2)*{^{t}},
 (0,1)*{\bullet}="a",
(8,1)*{\bu}="b",
\ar @{.>} "a";"b" <0pt>
\endxy}\Ea$, whose direction can be flipped up to the following sign factor,
$$
 \Ba{c}\resizebox{11mm}{!}{
\xy
(0,2.5)*{^{v_1}},
(8,2.5)*{^{v_2}},
(4,2)*{^{s}},
 (0,1)*{\bullet}="a",
(8,1)*{\bu}="b",
\ar @{.>} "a";"b" <0pt>
\endxy}\Ea
=(-1)^d
 \Ba{c}\resizebox{11mm}{!}{
\xy
(0,2.5)*{^{v_1}},
(8,2.5)*{^{v_2}},
(4,2)*{^{s}},
 (0,1)*{\bullet}="a",
(8,1)*{\bu}="b",
\ar @{<.} "a";"b" <0pt>
\endxy}\Ea
\ \ ,
\ \
 \Ba{c}\resizebox{11mm}{!}{
\xy
(0,2)*{^{v_1}},
(8,2)*{^{v_2}},
(4,1)*{^{t}},
 (0,0)*{\bullet}="a",
(8,0)*{\bu}="b",
\ar @{.>} "a";"b" <0pt>
\endxy}\Ea
=(-1)^d
 \Ba{c}\resizebox{11mm}{!}{
\xy
(0,2)*{^{v_1}},
(8,2)*{^{v_2}},
(4,1)*{^{t}},
 (0,0)*{\bullet}="a",
(8,0)*{\bu}="b",
\ar @{<.} "a";"b" <0pt>
\endxy}\Ea;
$$
\item[(iii)] wavy edges
  $\Ba{c}\resizebox{11mm}{!}{
\xy
(0,2)*{^{v_1}},
(8,2)*{^{v_2}},
 (0,0)*{\bullet}="a",
(8,0)*{\bu}="b",
\ar @{~>} "a";"b" <0pt>
\endxy}\Ea$   whose direction can be flipped up to the following sign factor,
$$
 \Ba{c}\resizebox{11mm}{!}{
\xy
(0,2)*{^{v_1}},
(8,2)*{^{v_2}},
 (0,0)*{\bullet}="a",
(8,0)*{\bu}="b",
\ar @{~>} "a";"b" <0pt>
\endxy}\Ea
=(-1)^d
 \Ba{c}\resizebox{11mm}{!}{
\xy
(0,2)*{^{v_1}},
(8,2)*{^{v_2}},
 (0,0)*{\bullet}="a",
(8,0)*{\bu}="b",
\ar @{<~} "a";"b" <0pt>
\endxy}\Ea
$$
\Ei
The cohomological degree of such a graph is given by
$$
|\Ga|=(d+1)(\# V(\Ga)-1) - d \# E_{sol}(\Ga)  + (1-d)\# E_{s\text{-}dot}(\Ga)  + (1-d)\# E_{t\text{-}dot}(\Ga) + (2-d)\# E_{wavy}(\Ga)
$$
where  $E_{sol}(\Ga)$, $E_{s\text{-}dot}(\Ga)$, $E_{t\text{-}dot}(\Ga)$ and $E_{wavy}(\Ga)$)
are, respectively, the sets of solid edges, dotted $s$-edges, $t$-edges and wavy edges of the graph $\Ga$.

\sip

A differential $\delta$ in the graded vector space $\rdGCD$  is defined by the formula
$$
\delta=\delta_\bu + \delta'
$$
where $\delta_\bu$ is the standard differential splitting the vertices with the help of a solid edge,
$$
\delta_\bu: \bu \rightsquigarrow \resizebox{9mm}{!}{
\xy
 (0,1)*{\bullet}="a",
(7,1)*{\bu}="b",
\ar @{->} "a";"b" <0pt>
\endxy}
$$
 while   $\delta'$ acts on edges of graphs $\Ga$  by changing their types as follows,
$$
\delta': \left\{
\Ba{l}
 \Ba{c}\resizebox{11mm}{!}{
\xy
(0,2.5)*{^{v_1}},
(8,2.5)*{^{v_2}},
 (0,1)*{\bullet}="a",
(8,1)*{\bu}="b",
\ar @{->} "a";"b" <0pt>
\endxy}\Ea
  \rightsquigarrow
 \Ba{c}\resizebox{11mm}{!}{
\xy
(0,2.5)*{^{v_1}},
(8,2.5)*{^{v_2}},
(4,2)*{^{t}},
 (0,1)*{\bullet}="a",
(8,1)*{\bu}="b",
\ar @{.>} "a";"b" <0pt>
\endxy}\Ea
\ - \
 \Ba{c}\resizebox{11mm}{!}{
\xy
(0,2.5)*{^{v_1}},
(8,2.5)*{^{v_2}},
(4,2)*{^{s}},
 (0,1)*{\bullet}="a",
(8,1)*{\bu}="b",
\ar @{.>} "a";"b" <0pt>
\endxy}\Ea
\vspace{1mm}
\\
\Ba{c}\resizebox{11mm}{!}{
\xy
(0,2.5)*{^{v_1}},
(8,2.5)*{^{v_2}},
(4,2)*{^{t}},
 (0,1)*{\bullet}="a",
(8,1)*{\bu}="b",
\ar @{.>} "a";"b" <0pt>
\endxy}\Ea
 \rightsquigarrow
\Ba{c}\resizebox{11mm}{!}{
\xy
(0,2.5)*{^{v_1}},
(8,2.5)*{^{v_2}},
 (0,1)*{\bullet}="a",
(8,1)*{\bu}="b",
\ar @{~>} "a";"b" <0pt>
\endxy}\Ea
\vspace{1mm}
\\
 \Ba{c}\resizebox{11mm}{!}{
\xy
(0,2.5)*{^{v_1}},
(8,2.5)*{^{v_2}},
(4,2)*{^{s}},
 (0,1)*{\bullet}="a",
(8,1)*{\bu}="b",
\ar @{.>} "a";"b" <0pt>
\endxy}\Ea
\rightsquigarrow
\Ba{c}\resizebox{11mm}{!}{
\xy
(0,2.5)*{^{v_1}},
(8,2.5)*{^{v_2}},
 (0,1)*{\bullet}="a",
(8,1)*{\bu}="b",
\ar @{~>} "a";"b" <0pt>
\endxy}\Ea
\Ea
\right.
$$

\subsubsection{{\bf Proposition} \cite{Z1}}\label{4: Prop on MZ's complex to dGC} {\em There is an explicit quasi-isomorphism of complexes
$$
\Ba{rccc}
z:& \rdGCD &\lon & \dGCD    \\
 & \Ga & \lon & z(\Ga)
\Ea
$$
where $z(\Ga)$ is obtained from $\Ga$ by changing its edges according to the rule
$$
z: \left\{
\Ba{l}
 \Ba{c}\resizebox{11mm}{!}{
\xy
(0,2.5)*{^{v_1}},
(8,2.5)*{^{v_2}},
 (0,1)*{\bullet}="a",
(8,1)*{\bu}="b",
\ar @{->} "a";"b" <0pt>
\endxy}\Ea
  \rightsquigarrow
   \Ba{c}\resizebox{11mm}{!}{
\xy
(0,2.5)*{^{v_1}},
(8,2.5)*{^{v_2}},
 (0,1)*{\bullet}="a",
(8,1)*{\bu}="b",
\ar @{->} "a";"b" <0pt>
\endxy}\Ea
\vspace{2mm}
\\
\Ba{c}\resizebox{11mm}{!}{
\xy
(0,2.5)*{^{v_1}},
(8,2.5)*{^{v_2}},
(4,2)*{^{t}},
 (0,1)*{\bullet}="a",
(8,1)*{\bu}="b",
\ar @{.>} "a";"b" <0pt>
\endxy}\Ea
 \rightsquigarrow
\Ba{c}\resizebox{13mm}{!}{\xy
(0,2)*{^{v_1}},
(7,5)*{^{v}},
(14,2)*{^{v_2}},
 (0,0)*{\bullet}="a",
 (7,3)*{\bu}="c",
(14,0)*{\bu}="b",
\ar @{->} "a";"c" <0pt>
\ar @{->} "b";"c" <0pt>
\endxy}
\Ea
\vspace{2mm}
\\
 \Ba{c}\resizebox{11mm}{!}{
\xy
(0,2.5)*{^{v_1}},
(8,2.5)*{^{v_2}},
(4,2)*{^{s}},
 (0,1)*{\bullet}="a",
(8,1)*{\bu}="b",
\ar @{.>} "a";"b" <0pt>
\endxy}\Ea
\rightsquigarrow
\Ba{c}\resizebox{13mm}{!}{
\xy
(0,2)*{^{v_1}},
(7,-1)*{^{v}},
(14,2)*{^{v_2}},
 (0,0)*{\bullet}="a",
 (7,-3)*{\bu}="c",
(14,0)*{\bu}="b",
\ar @{<-} "a";"c" <0pt>
\ar @{<-} "b";"c" <0pt>
\endxy}\Ea
\vspace{2mm}
\\
\Ba{c}\resizebox{11mm}{!}{
\xy
(0,2.5)*{^{v_1}},
(8,2.5)*{^{v_2}},
 (0,1)*{\bullet}="a",
(8,1)*{\bu}="b",
\ar @{~>} "a";"b" <0pt>
\endxy}\Ea
\rightsquigarrow
\Ba{c}\resizebox{18mm}{!}{\xy
(0,2)*{^{v_1}},
(7,5)*{^{v'}},
(14,2)*{^{v''}},
(21,5)*{^{v_2}},
 (0,0)*{\bullet}="a",
 (7,3)*{\bu}="c",
(14,0)*{\bu}="b",
(21,3)*{\bu}="d",
\ar @{->} "a";"c" <0pt>
\ar @{->} "b";"c" <0pt>
\ar @{->} "b";"d" <0pt>
\endxy}\Ea
+ (-1)^d
\Ba{c}\resizebox{18mm}{!}{\xy
(0,2)*{^{v_2}},
(7,5)*{^{v'}},
(14,2)*{^{v''}},
(21,5)*{^{v_1}},
 (0,0)*{\bullet}="a",
 (7,3)*{\bu}="c",
(14,0)*{\bu}="b",
(21,3)*{\bu}="d",
\ar @{->} "a";"c" <0pt>
\ar @{->} "b";"c" <0pt>
\ar @{->} "b";"d" <0pt>
\endxy}\Ea
\Ea
\right.
$$
}
\begin{proof}
This proposition was proven in \cite{Z1} in a much larger generality, in the case  of so called multioriented graph complexes. Let us show for completeness its short proof just for the case of ordinary directed graphs (the main idea is essentially the same as in \cite{Z1}).
The map $s$ in Proposition {\ref{4: Prop on 2-edge version of dGC}} factors through the map $z$ via the following morphism of complexes $$
\Ba{rccc}
s':& \rrdGCD &\lon & \rdGCD    \\
 & \Ga & \lon & s'(\Ga)
\Ea
$$
where $s'(\Ga)$ is obtained from $\Ga$ by changing its edges according to the rule
$$
s': \left\{
\Ba{l}
 \Ba{c}\resizebox{11mm}{!}{
\xy
(0,2.5)*{^{v_1}},
(8,2.5)*{^{v_2}},
 (0,1)*{\bullet}="a",
(8,1)*{\bu}="b",
\ar @{->} "a";"b" <0pt>
\endxy}\Ea
  \rightsquigarrow
   \Ba{c}\resizebox{11mm}{!}{
\xy
(0,2.5)*{^{v_1}},
(8,2.5)*{^{v_2}},
 (0,1)*{\bullet}="a",
(8,1)*{\bu}="b",
\ar @{->} "a";"b" <0pt>
\endxy}\Ea
\vspace{2mm}
\\
\Ba{c}\resizebox{11mm}{!}{
\xy
(0,2.5)*{^{v_1}},
(8,2.5)*{^{v_2}},
 (0,1)*{\bullet}="a",
(8,1)*{\bu}="b",
\ar @{.>} "a";"b" <0pt>
\endxy}\Ea
 \rightsquigarrow
 \Ba{c}\resizebox{11mm}{!}{
\xy
(0,2.5)*{^{v_1}},
(8,2.5)*{^{v_2}},
(4,2)*{^{t}},
 (0,1)*{\bullet}="a",
(8,1)*{\bu}="b",
\ar @{.>} "a";"b" <0pt>
\endxy}\Ea
-
\Ba{c}\resizebox{11mm}{!}{
\xy
(0,2.5)*{^{v_1}},
(8,2.5)*{^{v_2}},
(4,2)*{^{s}},
 (0,1)*{\bullet}="a",
(8,1)*{\bu}="b",
\ar @{.>} "a";"b" <0pt>
\endxy}\Ea
\Ea
\right.
$$
Choose a new basis in the complex $\rrdGCD$ in which dotted edges are replaced by the following
linear combinations,
$$
\Ba{c}\resizebox{11mm}{!}{
\xy
(0,2.5)*{^{v_1}},
(8,2.5)*{^{v_2}},
(4,2)*{^{(+)}},
 (0,1)*{\bullet}="a",
(8,1)*{\bu}="b",
\ar @{.>} "a";"b" <0pt>
\endxy}\Ea
:=
 \frac{1}{2}\left(
 \Ba{c}\resizebox{11mm}{!}{
\xy
(0,2.5)*{^{v_1}},
(8,2.5)*{^{v_2}},
(4,2)*{^{t}},
 (0,1)*{\bullet}="a",
(8,1)*{\bu}="b",
\ar @{.>} "a";"b" <0pt>
\endxy}\Ea
+
\Ba{c}\resizebox{11mm}{!}{
\xy
(0,2.5)*{^{v_1}},
(8,2.5)*{^{v_2}},
(4,2)*{^{s}},
 (0,1)*{\bullet}="a",
(8,1)*{\bu}="b",
\ar @{.>} "a";"b" <0pt>
\endxy}\Ea
\right),
\ \ \
\Ba{c}\resizebox{11mm}{!}{
\xy
(0,2.5)*{^{v_1}},
(8,2.5)*{^{v_2}},
(4,2)*{^{(-)}},
 (0,1)*{\bullet}="a",
(8,1)*{\bu}="b",
\ar @{.>} "a";"b" <0pt>
\endxy}\Ea
:=
 \frac{1}{2}\left(
 \Ba{c}\resizebox{11mm}{!}{
\xy
(0,2.5)*{^{v_1}},
(8,2.5)*{^{v_2}},
(4,2)*{^{t}},
 (0,1)*{\bullet}="a",
(8,1)*{\bu}="b",
\ar @{.>} "a";"b" <0pt>
\endxy}\Ea
-
\Ba{c}\resizebox{11mm}{!}{
\xy
(0,2.5)*{^{v_1}},
(8,2.5)*{^{v_2}},
(4,2)*{^{s}},
 (0,1)*{\bullet}="a",
(8,1)*{\bu}="b",
\ar @{.>} "a";"b" <0pt>
\endxy}\Ea
\right)
$$
and consider a filtration of the l.h.s. of the map $s'$ by the number of vertices plus the number of dotted edges (so that the induced differential on the associated graded is trivial), and the filtration of the r.h.s. by the number of vertices and the number of $(-)$-dotted edges so that the induced differential acts only on $(+)$-dotted edges by making them wavy,
$$
\delta_{ind}:
\Ba{c}\resizebox{11mm}{!}{
\xy
(0,2.5)*{^{v_1}},
(8,2.5)*{^{v_2}},
(4,2)*{^{(+)}},
 (0,1)*{\bullet}="a",
(8,1)*{\bu}="b",
\ar @{.>} "a";"b" <0pt>
\endxy}\Ea
 \rightsquigarrow
\Ba{c}\resizebox{11mm}{!}{
\xy
(0,2.5)*{^{v_1}},
(8,2.5)*{^{v_2}},
 (0,1)*{\bullet}="a",
(8,1)*{\bu}="b",
\ar @{~>} "a";"b" <0pt>
\endxy}\Ea
$$
At the next pages of the spectral sequences the map $s'$ becomes the isomorphism. Hence the map $s'$ is a quasi-isomorphism which in turn implies that $z$ is a quasi-isomorphism.
\end{proof}

\subsection{A further reduction of $\rdGCD$}
Let $\dGC_{d+1}^{\wedge\sim}$ be a subcomplex of $\rdGCD$ spanned by graphs having at least one
one dotted $t$-edge or at least wavy edge.

\subsubsection{\bf Proposition}\label{4: acyclicity of dGCD(wavy or wedge)}
{\em The subcomplex\, $\dGC_{d+1}^{\wedge\sim}$ is acyclic.}

\begin{proof}
Consider a filtration of $\dGC_{d+1}^{\wedge \sim}$ by the number of vertices. The induced differential on the associated graded complex $gr\dGC_{d+1}^{\wedge \sim}$ is $\delta'$ which acts only on edges by changing their types.
As the number of vertices in  $gr\dGC_{d+1}^{\wedge \sim}$ is preserved, we can consider its version  $gr^{marked}\dGC_{d+1}^{\wedge \sim}$ in which the vertices of graphs $\Ga$ are distinguished, say, an isomorphism $V(\Ga) \rar [\# \Ga]$ is fixed. By Maschke theorem, the acyclicity of the latter will imply the acyclicity of  $gr\dGC_{d+1}^{\wedge \sim}$ and hence
of $\dGC_{d+1}^{\wedge \sim}$.

\sip

Let $S$ be the set of graphs with labelled vertices of valencies $\geq 3$ whose edges have no types  assigned. We have an isomorphism of complexes
$$
gr^{marked}\dGC_{d+1}^{\wedge \sim}= \prod_{\Ga\in S} C_\Ga, \ \ \ \
C_\Ga:=\bigoplus_{e'\in \Ga}\left( C_{e'} \bigotimes_{e\in E(\Ga)\setminus e'} C_{e}\right)
$$
where the complexes $C_e$ and $C_{e'}$ control the graded vector spaces  all possible types which the underlying edges $e$ and $e'$ can have. Concretely,
$C_e$ is the reduced symmetric tensor algebra $\odot^{\bu \geq 1}C$ generated by
5-dimensional  complex
$$
 C=\text{span}\left\langle  \Ba{c}\resizebox{11mm}{!}{
 \xy
(0,2)*{^{v_1}},
(8,2)*{^{v_2}},
 (0,0)*{\bullet}="a",
(8,0)*{\bu}="b",
\ar @{->} "a";"b" <0pt>
\endxy}\Ea
\ , \
\Ba{c}\resizebox{11mm}{!}{
 \xy
(0,2)*{^{v_1}},
(8,2)*{^{v_2}},
 (0,0)*{\bullet}="a",
(8,0)*{\bu}="b",
\ar @{<-} "a";"b" <0pt>
\endxy}\Ea
\ , \
\Ba{c}\resizebox{11mm}{!}{
\xy
(0,2)*{^{v_1}},
(8,2)*{^{v_2}},
(4,1)*{^{s}},
 (0,0)*{\bullet}="a",
(8,0)*{\bu}="b",
\ar @{.>} "a";"b" <0pt>
\endxy}\Ea
 \ ,\
 \Ba{c}\resizebox{11mm}{!}{
\xy
(0,2)*{^{v_1}},
(8,2)*{^{v_2}},
(4,1)*{^{t}},
 (0,0)*{\bullet}="a",
(8,0)*{\bu}="b",
\ar @{.>} "a";"b" <0pt>
\endxy}\Ea
\ \ ,\ \
\Ba{c}\resizebox{11mm}{!}{
\xy
(0,2)*{^{v_1}},
(8,2)*{^{v_2}},
 (0,0)*{\bullet}="a",
(8,0)*{\bu}="b",
\ar @{~>} "a";"b" <0pt>
\endxy}\Ea
  \right\rangle
$$
whose cohomology group is generated as a graded symmetric tensor algebra by the 1-dimensional vector space with the basis vector given by the linear combination
$\Ba{c}\resizebox{10mm}{!}{
 \xy
(0,2.5)*{^{v_1}},
(8,2.5)*{^{v_2}},
 (0,1)*{\bullet}="a",
(8,1)*{\bu}="b",
\ar @{->} "a";"b" <0pt>
\endxy}\Ea\hspace{-1mm}
- (-1)^d\hspace{-1mm}
\Ba{c}\resizebox{10mm}{!}{
 \xy
(0,2.5)*{^{v_1}},
(8,2.5)*{^{v_2}},
 (0,1)*{\bullet}="a",
(8,1)*{\bu}="b",
\ar @{<-} "a";"b" <0pt>
\endxy}\Ea$. One the other hand, one has $C_{e'}\simeq C'\ot \odot^{\bu \geq 0}C$, where $C'$ is
a  two-dimensional complex
$$
C'=\text{span}\left\langle  \Ba{c}\resizebox{11mm}{!}{
\xy
(0,2)*{^{v_1}},
(8,2)*{^{v_2}},
(4,1)*{^{t}},
 (0,0)*{\bullet}="a",
(8,0)*{\bu}="b",
\ar @{.>} "a";"b" <0pt>
\endxy}\Ea
\ \ ,\ \
\Ba{c}\resizebox{11mm}{!}{
\xy
(0,2)*{^{v_1}},
(8,2)*{^{v_2}},
 (0,0)*{\bullet}="a",
(8,0)*{\bu}="b",
\ar @{~>} "a";"b" <0pt>
\endxy}\Ea
  \right\rangle.
$$
As $C'$ is acyclic, each complex $C_\Ga$ is acyclic implying in turn the acyclicity
of $gr^{marked}\dGC_{d+1}^{\wedge \sim}$. The proposition is proven.
\end{proof}

The  quotient of $\rdGCD$ by the above acyclic subcomplex is isomorphic to the reduced directed complex $\rrdGCD$ introduced earlier, i.e.\ the two reduced models of $\dGCD$ are related
by a short exact sequence
$$
0\lon \dGC_{d+1}^{\wedge \sim} \lon  \rdGCD \stackrel{\tilde{\pi}}{\lon} \rrdGCD \lon 0
$$
where $\tilde{\pi}$
is a quasi-isomorphism.

\sip

 Consider next a subcomplex $\dGC_{d+1}^{\wedge}$  of $\dGCD$ spanned by graphs having at least one bivalent target.

\subsubsection{\bf Proposition}\label{4: Prop on acyclicity of GCvee} {\em
The complex $\dGC_{d+1}^{\wedge}$ is acyclic.}

\begin{proof} There is an explicit morphism from the reduced directed graph complex $\rrdGCD$
into the following quotient complex
$$
g: \rrdGCD \lon \dGCD/\dGC_{d+1}^{\wedge}
$$
which preserves vertices and solid edges of the generating graphs
$\Ga\in \rrdGCD$ while sending every dotted edge of $\Ga$ into a bivalent source
$$
g:\
\xy
(0,2)*{^{v_1}},
(8,2)*{^{v_2}},
 (0,0)*{\bullet}="a",
(8,0)*{\bu}="b",
\ar @{.>} "a";"b" <0pt>
\endxy
\lon
-
\xy
(0,2)*{^{v_1}},
(8,2)*{^{v}},
(16,2)*{^{v_2}},
 (0,0)*{\bullet}="a",
 (8,0)*{\bu}="c",
(16,0)*{\bu}="b",
\ar @{<-} "a";"c" <0pt>
\ar @{<-} "b";"c" <0pt>
\endxy,
$$
where the new solid edges on the r.h.s\ are ordered from left to right, while the new vertex
$v$ is assumed to be adjoined into the given ordering of vertices of $\Ga$ as the last one $v_1\wedge v_2\wedge v$.
This morphism is obviously an isomorphism implying the claim.
\end{proof}

Similarly one proves acyclicity of the subcomplex $\GC_{d+1}^{\vee}\hook \dGC_{d+1}$ spanned by graphs having at least one
bivalent source.

\sip

Let  $\mathsf{GC}_{d+1}^{\vee + \wedge}$ (resp.\, $\GCD^{\vee \cdot \wedge}$)    be a subcomplex of $\dGCD$ generated by graphs with at least one bivalent source {\em or}\, (resp., {\em and}) bivalent target; we shall see such subcomplexes in applications below. There is a short exact sequence of complexes
$$
0\lon \mathsf{GC}_{d+1}^{\wedge\cdot \vee} \lon \mathsf{GC}_{d+1}^{\wedge} \oplus
\mathsf{GC}_{d+1}^{\vee} \lon \mathsf{GC}_{d+1}^{\vee + \wedge} \lon 0,
$$
As the middle term is acyclic, we obtain the following

\subsubsection{\bf Corollary}
$
H^\bu(\mathsf{GC}_{d+1}^{\vee + \wedge})=H^{\bu+1}(\mathsf{GC}_{d+1}^{\wedge\cdot \vee})$.

\sip


\subsubsection{Proof of Theorem {\ref{4: Theorem rTGCd versis TGCd}}}  It is enough to show
$H^\bu(\rTGCd)=H^\bu(\TGCd)$. Indeed, the targeted complex $\TGCd$ contains an acyclic complex
$\GC_{d+1}^\wedge$ so that the projection $\pi$ in the following
 short exact sequence
$$
0\lon \GC_{d+1}^\wedge  \lon \TGCd \stackrel{\pi}{\lon} \rTGCd \lon 0.
$$
is a quasi-isomorphism. The quotient complex is obviously isomorphic to the complex  $\rTGCd$ introduced in \S 4.1 (hence the notation). The theorem is proven.


\mip

{\large
\section{\bf Directed graph complex: its auxiliary extensions}
}

\subsection{An auxiliary extension of $\dGC_{d+1}$ with two types of vertices and one type of edges}
Let
$\ddGCdd$ be a graded vector space generated by directed connected graphs $\Ga$ such that
\Bi
\item[(i)] $\Ga$ has two types of vertices, white and black ones, and only one type of edges --- the solid edges with a fixed direction; white vertices  are allowed to have only incoming edges, i.e. \ they are all targets;
\item[(ii)] black and white vertices of $\Ga$ are at least bivalent, and at least one (black or white) vertex is of valency $\geq 3$, no passing black vertices are allowed, e.g.
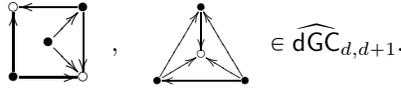
 $$
 \Ba{c}\resizebox{12mm}{!}{
\xy
 (0,0)*{\bullet}="a",
(0,10)*{\circ}="b",
(10,0)*{\circ}="c",
(10,10)*{\bu}="d",
(5,5)*{\bu}="o",
\ar @{->} "a";"b" <0pt>
\ar @{->} "a";"c" <0pt>
\ar @{->} "d";"c" <0pt>
\ar @{->} "d";"b" <0pt>
\ar @{->} "o";"d" <0pt>
\ar @{->} "o";"c" <0pt>
\endxy}
\Ea
\ , \
\Ba{c}\resizebox{14mm}{!}{
\xy
 (0,0)*{\circ}="a",
(0,8)*{\bullet}="b",
(-7.5,-4.5)*{\bu}="c",
(7.5,-4.5)*{\bu}="d",
\ar @{<-} "a";"b" <0pt>
\ar @{<-} "a";"c" <0pt>
\ar @{<-} "b";"c" <0pt>
\ar @{->} "d";"c" <0pt>
\ar @{<-} "b";"d" <0pt>
\ar @{->} "d";"a" <0pt>
\endxy}
\Ea
\in \ddGCdd.
$$

\item[(iii)]  The cohomological degree of $\Ga\in \ddGCdd$ is given by
$$
|\Ga|=(d+1)(\# V_{\bu}(\Ga) -1) + d\#V_{\circ}(\Ga) + (-d) \# E(\Ga).
$$
where  $V_{\bu}(\Ga)$ (resp., $V_{\circ}(\Ga)$) is the set of black (resp., white) vertices of $\Ga$ and $E(\Ga)$ is the set of edges of $\Ga$.
\Ei

\sip

The vector space $\ddGCdd$ contains $\dGCD$ as a subspace generated by graphs with all vertices black.

\sip

Consider a degree one linear map
\Beq\label{4: d in ddGCdd}
\Ba{rccl}
d: & \ddGCdd &\lon & \ddGCdd\\
        &\Ga &\lon &d \Ga :=\underbrace{\sum_{v\in V_{\bu}(\Ga)} \delta_{\bu,v} \Ga}_{\delta_{\bu}\Ga} +
\underbrace{\sum_{v\in V_{\circ}(\Ga)} d_{\circ,v} \Ga}_{d_{\circ}\Ga},
\Ea
\Eeq
where
\Bi
\item[(i)] the operator $\delta_{\bu,v}$ acts on the black vertex $v$ by splitting it into a pair of new black vertices,
$$
\delta_{\bu,v}:\bu \rightsquigarrow \resizebox{10mm}{!}{
\xy
 (0,1)*{\bullet}="a",
(7,1)*{\bu}="b",
\ar @{->} "a";"b" <0pt>
\endxy}
$$
and redistributing the edges attached to $v$ among the new vertices in all possible ways
(as in (\ref{2: delta_v splits vertex into two})) such that no new univalent vertices are created;

\item[(ii)]
the operator $d_{\circ,v}$ acts on the white vertex $v$ of $\Ga$ by substituting into $v$ the following linear combination of graphs
$$
d_{\circ,v}: \circ \rightsquigarrow
\ \bu
+ \sum_{k=1}^\infty \frac{1}{k!}
\resizebox{11mm}{!}{
\xy
 (-5,-1)*{\circ}="a1",
  (5,-1)*{\circ}="a2",
   (-2,-1)*{\circ}="a3",
    (2,-1)*{...},
   (0,8)*{\bu}="b",
(0,-5)*{\underbrace{\hspace{12mm}}_k},
\ar @{<-} "a1";"b" <0pt>
\ar @{<-} "a2";"b" <0pt>
\ar @{<-} "a3";"b" <0pt>
\endxy}
$$
and redistributing the edges attached to $v$ among the vertices of each summand on the right hand side in such a way that the new vertices are at least bivalent; notice that this procedure can not create white vertices which are not targets. The operation $\delta_{\circ,v}\Ga$ produces always a finite linear combination  of graphs.
\Ei

For each $\Ga\in \ddGCdd$ the orientation of $\delta \Ga$ is defined in a full analogy to the cases of $\GCd$ and $\dGCD$.

\subsubsection{\bf Lemma}\label{4: Lemma on d^2=0}  $d^2=0$.

\begin{proof} Let $\mathsf{fc}\ddGCdd$ be a (full connected)  extension of $\ddGCdd$ in which we allow generating  graphs to have vertices of any valency, including 0-valent and 1-valent ones. We define an operator $\delta$ in $\mathsf{fc}\ddGCdd$ by formally the same formulae as above (i.e.\ via substitutions of appropriate linear combinations of graphs into vertices) except that  the edges gets redistributed among newly created vertices in {\em all}\, possible ways, i.e.\ we do not impose any vanishing conditions (such as no new 1-valent vertices etc).
The lemma is essentially proven once we show that $\delta^2=0$ in $\mathsf{fc}\ddGCdd$ which in turn is proven once we show that the operation $\delta^2$ is zero when applied to the one-vertex graphs,
$$
\delta^2(\bu)=0, \ \ \ \delta^2(\circ)=0.
$$
The first equality is obvious (as it does not involve the new white vertices). Let us check the second one,
\Beqrn
\delta^2(\circ) &=& \delta_\bu\left(  \sum_{k=0}^\infty \frac{1}{k!}
\resizebox{11mm}{!}{
\xy
 (-5,-1)*{\circ}="a1",
  (5,-1)*{\circ}="a2",
   (-2,-1)*{\circ}="a3",
    (2,-1)*{...},
   (0,8)*{\bu}="b",
(0,-5)*{\underbrace{\hspace{12mm}}_k},
\ar @{<-} "a1";"b" <0pt>
\ar @{<-} "a2";"b" <0pt>
\ar @{<-} "a3";"b" <0pt>
\endxy}  \right)
+
\delta_\circ\left(  \sum_{p=0}^\infty \frac{1}{k!}
\resizebox{11mm}{!}{
\xy
 (-5,-1)*{\circ}="a1",
  (5,-1)*{\circ}="a2",
   (-2,-1)*{\circ}="a3",
    (2,-1)*{...},
   (0,8)*{\bu}="b",
(0,-5)*{\underbrace{\hspace{12mm}}_p},
\ar @{<-} "a1";"b" <0pt>
\ar @{<-} "a2";"b" <0pt>
\ar @{<-} "a3";"b" <0pt>
\endxy}  \right)
\\
&=&\left(
\sum_{k=p+q\atop p\geq 0, q\geq 0} \frac{1}{k!}\frac{k!}{p!q!}
\resizebox{20mm}{!}{
\xy
 (7,-1)*{\circ}="a1'",
  (15,-1)*{\circ}="a2'",
   (10,-1)*{\circ}="a3'",
    (13,-1)*{...},
   (10,8)*{\bu}="b'",
(11,-5)*{\underbrace{\hspace{9mm}}_q},
 (-5,-1)*{\circ}="a1",
  (3,-1)*{\circ}="a2",
   (-2,-1)*{\circ}="a3",
    (1,-1)*{...},
   (0,8)*{\bu}="b",
(-1,-5)*{\underbrace{\hspace{9mm}}_p},
\ar @{<-} "a1";"b" <0pt>
\ar @{<-} "a2";"b" <0pt>
\ar @{<-} "a3";"b" <0pt>
\ar @{<-} "a1'";"b'" <0pt>
\ar @{<-} "a2'";"b'" <0pt>
\ar @{<-} "a3'";"b'" <0pt>
\ar @{<-} "b'";"b" <0pt>
\endxy}  \right)
-
\left(\sum_{p\geq 0, q\geq 0} \frac{1}{p!q!}
\resizebox{20mm}{!}{
\xy
 (7,-1)*{\circ}="a1'",
  (15,-1)*{\circ}="a2'",
   (10,-1)*{\circ}="a3'",
    (13,-1)*{...},
   (10,8)*{\bu}="b'",
(11,-5)*{\underbrace{\hspace{9mm}}_q},
 (-5,-1)*{\circ}="a1",
  (3,-1)*{\circ}="a2",
   (-2,-1)*{\circ}="a3",
    (1,-1)*{...},
   (0,8)*{\bu}="b",
(-1,-5)*{\underbrace{\hspace{9mm}}_p},
\ar @{<-} "a1";"b" <0pt>
\ar @{<-} "a2";"b" <0pt>
\ar @{<-} "a3";"b" <0pt>
\ar @{<-} "a1'";"b'" <0pt>
\ar @{<-} "a2'";"b'" <0pt>
\ar @{<-} "a3'";"b'" <0pt>
\ar @{<-} "b'";"b" <0pt>
\endxy}  \right)=0.
\Eeqrn

Next we notice that the subspace of $\mathsf{fc}\ddGCdd$ spanned by graphs having at least one vertex of valency $\leq 1$   is a subcomplex; the associated quotient complex is precisely $\ddGCdd$.
The Lemma is proven.
\end{proof}

\subsection{ An extension of $\dGCD$ with two types of vertices and two types of edges}\label{4: proof of delta^2} Consider for a moment a further extension of $\mathsf{fc}\ddGCdd$ :
let $\mathsf{fc}\dGC_{d,d+1}^{+}$ be a graded vector space generated by graphs with {\em two}\, types of vertices  exactly as in the case of $\mathsf{fc}\ddGCdd$  and with {\em two}\, types of edges, solid edges $
 \xy
 (0,0)*{\bullet}="a",
(7,0)*{\bu}="b",
\ar @{->} "a";"b" <0pt>
\endxy$ whose direction is fixed, and dotted ones
 $\xy
 (0,0)*{\bullet}="a",
(7,0)*{\bu}="b",
\ar @{.>} "a";"b" <0pt>
\endxy
$
whose direction can be flipped as in (\ref{2: symmetry of dotted edges}),
e.g.
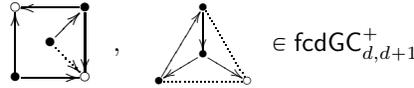
 $$
  \Ba{c}\resizebox{12mm}{!}{
\xy
 (0,0)*{\bullet}="a",
(0,10)*{\circ}="b",
(10,0)*{\circ}="c",
(10,10)*{\bu}="d",
(5,5)*{\bu}="o",
\ar @{->} "a";"b" <0pt>
\ar @{->} "a";"c" <0pt>
\ar @{->} "d";"c" <0pt>
\ar @{->} "d";"b" <0pt>
\ar @{->} "o";"d" <0pt>
\ar @{.>} "o";"c" <0pt>
\endxy}
\Ea
\ , \
 \Ba{c}\resizebox{14mm}{!}{
\xy
 (0,0)*{\bullet}="a",
(0,8)*{\bullet}="b",
(-7.5,-4.5)*{\bu}="c",
(7.5,-4.5)*{\circ}="d",
\ar @{<-} "a";"b" <0pt>
\ar @{->} "a";"c" <0pt>
\ar @{<-} "b";"c" <0pt>
\ar @{.} "d";"c" <0pt>
\ar @{.} "b";"d" <0pt>
\ar @{<-} "d";"a" <0pt>
\endxy}
\Ea\in \mathsf{fc}\dGC_{d,d+1}^{+}
 $$

The cohomological degree of $\Ga\in \mathsf{fc}\dGC_{d,d+1}^{+}$ is given by
\Beq\label{3: deg of graphs in OGCdd}
|\Ga|=(d+1)(\# V_{\bu}(\Ga) -1) + d\#V_{\circ}(\Ga) -d \# E_{sol}(\Ga) +  (1-d) \# E_{dot}(\Ga).
\Eeq
where  $V_{\bu}(\Ga)$ (resp., $V_{\circ}(\Ga)$) is the set of black (resp., white) vertices of $\Ga$ and $E_{sol}(\Ga)$ (resp. $E_{dot}(\Ga)$)  is the set of solid (resp., dotted)  edges of $\Ga$.

\sip

 Consider a degree 1 endomorphism of $\mathsf{fc}\dGC_{d,d+1}^{+}$ given by
$$
\Ba{rccl}
\p: & \mathsf{fc}\dGC_{d,d+1}^{+} &\lon & \mathsf{fc}\dGC_{d,d+1}^{+}\\
        &\Ga &\lon &
         \Ga :=(\delta_{\bu} + \delta_{\circ})\Ga +  \delta_{\circ\circ}\Ga +\delta'\Ga
\Ea
$$
where $d:=\delta_\bu  + \delta_\circ$ is defined by exactly the same formulae as in the case of $\mathsf{fc}\ddGCdd$ discussed just above, the operator $\delta_{\circ\circ}$ acts on white vertices only by splitting them as follows
$$
\delta_{\circ\circ}\Ga:=\sum_{v\in V_{\circ}(\Ga)} \delta_{\circ\circ,v} \Ga, \ \ \
 \delta_{\circ\circ,v}: \circ \rightsquigarrow  \frac{1}{2}
\resizebox{11mm}{!}{
\xy
 (0,1)*{\circ}="a",
(7,1)*{\circ}="b",
\ar @{.>} "a";"b" <0pt>
\endxy},
$$
and the operator $\delta'$ acts on solid edges by making them dotted (as usual),
  $$
\delta'\Ga:= \sum_{e\in E_{sol}(\Ga)} \delta'_e \Ga \ \ \text{where}\ \
\delta'_e:   \resizebox{10mm}{!}{
\xy
%
 (0,2)*{}="a",
(7,2)*{}="b",
\ar @{->} "a";"b" <0pt>
\endxy}
\rightsquigarrow
 \resizebox{10mm}{!}{
\xy
%
 (0,2)*{}="a",
(7,2)*{}="b",
\ar @{.>} "a";"b" <0pt>
\endxy}\, .
$$

  \subsubsection{\bf Lemma}\label{5: Lemma p^2=0}  $\p^2=0$.

\begin{proof}
 As $d^2=\delta_{\circ\circ}^2=[\delta', \delta_{\circ\circ}]=0$, the Lemma is proven once we show that
 $[\delta_{\circ\circ}+ \delta',d]=0$.  It is enough to check the latter condition by applying its left hand side to the one white vertex graph. We have
\Beqrn
(\delta'\cdot d + \delta_{\circ\circ}\cdot d + d\cdot \delta_{\circ\circ})\,\circ &=&
(\delta' + \delta_{\circ\circ})\left(  \sum_{k=0}^\infty \frac{1}{k!}
\resizebox{11mm}{!}{
\xy
 (-5,-1)*{\circ}="a1",
  (5,-1)*{\circ}="a2",
   (-2,-1)*{\circ}="a3",
    (2,-1)*{...},
   (0,8)*{\bu}="b",
(0,-5)*{\underbrace{\hspace{12mm}}_k},
\ar @{<-} "a1";"b" <0pt>
\ar @{<-} "a2";"b" <0pt>
\ar @{<-} "a3";"b" <0pt>
\endxy}   \right)
 + d\left(\frac{1}{2}
\resizebox{10mm}{!}{
\xy
 (0,1)*{\circ}="a",
(7,1)*{\circ}="b",
\ar @{.>} "a";"b" <0pt>
\endxy}\right)\\
&=& \left(
 \sum_{k=1}^\infty \frac{1}{k!}
\resizebox{11mm}{!}{
\xy
 (-5,-1)*{\circ}="a1",
  (5,-1)*{\circ}="a2",
   (-2,-1)*{\circ}="a3",
    (2,-1)*{...},
   (0,8)*{\bu}="b",
(0,-5)*{\underbrace{\hspace{12mm}}_k},
\ar @{.} "a1";"b" <0pt>
\ar @{<-} "a2";"b" <0pt>
\ar @{<-} "a3";"b" <0pt>
\endxy}
+
 \sum_{k=1}^\infty \frac{1}{k!}
\resizebox{15mm}{!}{
\xy
(-10,-1)*{\circ}="a0",
 (-5,-1)*{\circ}="a1",
  (5,-1)*{\circ}="a2",
   (-2,-1)*{\circ}="a3",
    (2,-1)*{...},
   (0,8)*{\bu}="b",
(0,-5)*{\underbrace{\hspace{12mm}}_k},
\ar @{.} "a1";"a0" <0pt>
\ar @{<-} "a1";"b" <0pt>
\ar @{<-} "a2";"b" <0pt>
\ar @{<-} "a3";"b" <0pt>
\endxy} \right) -
\left(
 \sum_{k=1}^\infty \frac{1}{k!}
\resizebox{11mm}{!}{
\xy
 (-5,8)*{\circ}="a1",
  (5,-1)*{\circ}="a2",
   (-2,-1)*{\circ}="a3",
    (2,-1)*{...},
   (0,8)*{\bu}="b",
(2,-5)*{\underbrace{\hspace{9mm}}_k},
\ar @{.} "a1";"b" <0pt>
\ar @{<-} "a2";"b" <0pt>
\ar @{<-} "a3";"b" <0pt>
\endxy}
+
 \sum_{k=1}^\infty \frac{1}{k!}
\resizebox{12mm}{!}{
\xy
(-5,8)*{\circ}="a0",
 (-5,-1)*{\circ}="a1",
  (5,-1)*{\circ}="a2",
   (-2,-1)*{\circ}="a3",
    (2,-1)*{...},
   (0,8)*{\bu}="b",
(0,-5)*{\underbrace{\hspace{12mm}}_k},
\ar @{.} "a1";"a0" <0pt>
\ar @{<-} "a1";"b" <0pt>
\ar @{<-} "a2";"b" <0pt>
\ar @{<-} "a3";"b" <0pt>
\endxy} \right)\\
&=&0.
\Eeqrn
The claim is proven.
\end{proof}

\sip

The complex $ \mathsf{fc}\dGC_{d,d+1}^{+}$ contains a subcomplex $\sY^{\leq 2}$ generated by
graphs having at least one vertex of valency $\leq 2$. We denote by
 $\dddGCdd$ the associated  quotient complex,
 $\mathsf{fc}\dGC_{d,d+1}^{+}/\sY^{\leq 2}$, which is
 generated by {\em at least  trivalent}\, graphs with two types of vertices and two types of edges. We shall use that quotient complex in \S 6 below.

\mip


{\large
\section{\bf Graph complexes interpolating between $\GCd$ and $\OOGCd$/ $\rTGCd$}
}

\sip

\subsection{Reduced versions of $\TGCd$ and $\OGCd$}  We have described above the reduced models $\rOGCd$ and $\rTGCd$ of the oriented and targeted graph complexes (see \S 2.5 and \S 4.1 respectively). The arguments and constructions given below work  well  for both these complexes. Thus we use from now on the notation $\lGC$ where the parameter $\la$ takes values
in the set of two abbreviations, $\{\mathsf{t}, \mathsf{or}\}$.

\subsection{A complex with two types of vertices} Let $\llGCdd$  be a subcomplex of $\dddGCdd$ spanned, for $\lambda=\mathsf{t}$
by graphs having at least one vertex (of any type) with no outgoing solid edges, and,
for $\lambda=\mathsf{or}$ by graphs which have no closed directed paths of solid edges.

\sip

For reader's convenience, let us summarize briefly the main facts about $\llGCdd$. This is a
 graded vector space generated by connected graphs $\Ga$ with
\Bi
\item[(i)]  two types of vertices, white and black ones, which have valency $\geq 3$;
\item[(ii)]  two types of edges: dotted ones whose direction can be flipped and multiplied by $(-1)^d$,
and solid ones which have a fixed direction; white vertices  are not allowed to have outgoing solid edges attached;
\item[(iii)] for $\la=\mathsf{t}$ every graph $\Ga$ has at least one target, that is a (black or white) vertex with no outgoing solid edges; for $\la=\mathsf{or}$ every $\Ga$ has no closed directed paths of solid edges (``solid wheels").
As direction on dotted edges can be flipped, we often do not show such directions in our pictures, for example,
 $$
 \Ba{c}\resizebox{14mm}{!}{
\xy
 (0,0)*{\bullet}="a",
(0,8)*{\circ}="b",
(-7.5,-4.5)*{\circ}="c",
(7.5,-4.5)*{\circ}="d",
\ar @{.} "a";"b" <0pt>
\ar @{.} "a";"c" <0pt>
\ar @{.} "b";"c" <0pt>
\ar @{.} "d";"c" <0pt>
\ar @{.} "b";"d" <0pt>
\ar @{.} "d";"a" <0pt>
\endxy}
\Ea
\ , \
\Ba{c}\resizebox{14mm}{!}{
\xy
 (0,0)*{\bullet}="a",
(0,8)*{\bullet}="b",
(-7.5,-4.5)*{\circ}="c",
(7.5,-4.5)*{\circ}="d",
\ar @{<-} "a";"b" <0pt>
\ar @{->} "a";"c" <0pt>
\ar @{->} "b";"c" <0pt>
\ar @{.} "d";"c" <0pt>
\ar @{.} "b";"d" <0pt>
\ar @{<-} "d";"a" <0pt>
\endxy}
\Ea
 \ , \
  \Ba{c}\resizebox{14mm}{!}{
\xy
 (0,8)*{\bullet}="a",
(-7.5,-4.5)*{\bullet}="c",
(7.5,-4.5)*{\circ}="d",
\ar@/^-0.8pc/@{.}"a";"c" <0pt>
\ar@/^0.8pc/@{.}"a";"d" <0pt>
\ar @{<-} "a";"c" <0pt>
\ar @{.} "d";"c" <0pt>
\ar @{->} "a";"d" <0pt>
\endxy}
\Ea\in \OOGCdd\subset \TTGCdd, \  \Ba{c}\resizebox{14mm}{!}{
\xy
 (0,0)*{\bullet}="a",
(0,8)*{\bullet}="b",
(-7.5,-4.5)*{\bu}="c",
(7.5,-4.5)*{\circ}="d",
\ar @{<-} "a";"b" <0pt>
\ar @{->} "a";"c" <0pt>
\ar @{<-} "b";"c" <0pt>
\ar @{.} "d";"c" <0pt>
\ar @{.} "b";"d" <0pt>
\ar @{<-} "d";"a" <0pt>
\endxy}
\Ea\in  \TTGCdd\setminus \OOGCdd.
$$

\item[(iv)]  The cohomological degree of $\Ga\in \llGCdd$ is given by
\Beq\label{3: deg of graphs in llGCdd}
|\Ga|=(d+1)(\# V_{\bu}(\Ga) -1) + d\#V_{\circ}(\Ga) + (1-d)\# E_{dot}(\Ga) + (-d) \# E_{sol}(\Ga).
\Eeq
where  $V_{\bu}(\Ga)$ (resp., $V_{\circ}(\Ga)$) is the set of black (resp., white) vertices of $\Ga$ and $E_{sol}(\Ga)$ (resp., $E_{dot}(\Ga)$) is the set of solid (resp., dotted) edges of $\Ga$.
\Ei
The differential in $\llGCdd$ does depend on the particular value of the parameter $\lambda$ and is defined as linear map
\Beq\label{3: delta in GCd,d+1}
\Ba{rccl}
\delta: & \llGCdd &\lon & \llGCdd\\
        &\Ga &\lon &\delta \Ga :=\underbrace{\sum_{v\in V_{\bu}(\Ga)} \delta_{\bu,v} \Ga}_{\delta_{\bu}\Ga} +
\underbrace{\sum_{v\in V_{\circ}(\Ga)} \delta_{\circ,v} \Ga}_{\delta_{\circ}\Ga}
+
\underbrace{\sum_{e\in E_{sol}(\Ga)} \delta_e \Ga}_{\delta' \Ga},
\Ea
\Eeq
where
\Bi
\item[(i)] the operator $\delta_{\bu,v}$ acts on the black vertex $v$ by splitting it into a pair of new black vertices,
$$
\delta_{\bu,v}:\bu \rightsquigarrow \resizebox{10mm}{!}{
\xy
 (0,1)*{\bullet}="a",
(7,1)*{\bu}="b",
\ar @{->} "a";"b" <0pt>
\endxy}
$$
and redistributing the edges attached to $v$ among the new vertices in all possible ways
(as in (\ref{2: delta_v splits vertex into two})) while keeping the new vertices at least trivalent; thus $\delta_{\bu, v}$ acts trivially on trivalent black vertices.

\item[(ii)]
the operator $\delta_{\circ,v}$ acts on a white vertex $v$ of $\Ga$ by replacing $v$
with the following linear combination of graphs
$$
\delta_{\circ,v}: \circ \rightsquigarrow  \frac{1}{2}
\resizebox{11mm}{!}{
\xy
 (0,1)*{\circ}="a",
(7,1)*{\circ}="b",
\ar @{.>} "a";"b" <0pt>
\endxy}
\ +
\ \bu
+ \sum_{k=1}^\infty \frac{1}{k!}
\resizebox{11mm}{!}{
\xy
 (-5,-1)*{\circ}="a1",
  (5,-1)*{\circ}="a2",
   (-2,-1)*{\circ}="a3",
    (2,-1)*{...},
   (0,8)*{\bu}="b",
(0,-5)*{\underbrace{\hspace{12mm}}_k},
\ar @{<-} "a1";"b" <0pt>
\ar @{<-} "a2";"b" <0pt>
\ar @{<-} "a3";"b" <0pt>
\endxy}
$$
and redistributing the edges attached to $v$ among the vertices of each summand in such a way that thew new vertices are at least trivalent; in particular,
$\delta_{\circ,v}\Ga=0$ if $v\in V_{\circ}(\Ga)$ is trivalent.

\item[(iii)] the operator  $\delta'_e$ acts on a solid edges of $\Ga$ by changing its type (as in the case of the complex $\lGC$)
$$
\delta_e:
\resizebox{11mm}{!}{
\xy
(0,3)*{},
(8,3)*{},
 (0,1)*{\bullet}="a",
(8,1)*{\circ}="b",
\ar @{->} "a";"b" <0pt>
\endxy}
\lo
\resizebox{11mm}{!}{
\xy
(0,3)*{},
(8,3)*{},
 (0,1)*{\bullet}="a",
(8,1)*{\circ}="b",
\ar @{.>} "a";"b" <0pt>
\endxy}\ \ \text{or}
\ \
\resizebox{11mm}{!}{
\xy
%
 (0,1)*{\bullet}="a",
(8,1)*{\bu}="b",
\ar @{->} "a";"b" <0pt>
\endxy}
\lo
\resizebox{11mm}{!}{
\xy
(0,3)*{},
(8,3)*{},
 (0,1)*{\bullet}="a",
(8,1)*{\bu}="b",
\ar @{.>} "a";"b" <0pt>
\endxy},
$$
\Ei

For each $\Ga$ the orientation of $\delta \Ga$ is defined in a full analogy to the
cases of of $\GCd$ and $\OGCd$. Fur future use we split the operator $\delta_0$
into two summands
$$
\delta_\circ = \delta_{\circ\circ} + \deltabw
$$
where $\delta_{\circ\circ}$ is the part which does not create new black vertex and
$$
{\deltabw}\Ga:= \sum_{v\in V_{\circ}(\Ga)}
\delta_{\circ\hspace{-0.5mm}-\hspace{-1mm}\bu,v} \Ga \ \ \
\text{with}\ \ \ \delta_{\circ\hspace{-0.5mm}-\hspace{-1mm}\bu,v}:
\circ \rightsquigarrow   \bu
+ \sum_{k=1}^\infty \frac{1}{k!}
\resizebox{11mm}{!}{
\xy
 (-5,-1)*{\circ}="a1",
  (5,-1)*{\circ}="a2",
   (-2,-1)*{\circ}="a3",
    (2,-1)*{...},
   (0,8)*{\bu}="b",
(0,-5)*{\underbrace{\hspace{12mm}}_k},
\ar @{<-} "a1";"b" <0pt>
\ar @{<-} "a2";"b" <0pt>
\ar @{<-} "a3";"b" <0pt>
\endxy}=
\sum_{k=0}^\infty \frac{1}{k!}
\resizebox{11mm}{!}{
\xy
 (-5,-1)*{\circ}="a1",
  (5,-1)*{\circ}="a2",
   (-2,-1)*{\circ}="a3",
    (2,-1)*{...},
   (0,8)*{\bu}="b",
(0,-5)*{\underbrace{\hspace{12mm}}_k},
\ar @{<-} "a1";"b" <0pt>
\ar @{<-} "a2";"b" <0pt>
\ar @{<-} "a3";"b" <0pt>
\endxy}
$$

\subsubsection{\bf Acyclicity Lemma}\label{6: Lemma on acyclicity of llGCdd}
{\em The complex $\llGCdd$ is acyclic.}

\begin{proof}
 Consider a filtration of $\llGCdd$
 by
the total number of vertices plus the
number of dotted edges. The induced differential in the associated graded complex $gr\llGCdd$  just makes white vertices (if any) into black ones,
$\circ \rightsquigarrow \bu$.  One can consider a version $\GC^\la_{marked}$ of $gr\llGCdd$ in which
all edges and vertices are marked but the {\em type}\,  of vertices is not fixed, i.e they can in general be decorated in black or white colours, and the differential acts only on decorations. By Maschke Theorem, the acyclicity of $\GC_{marked}^\la$ implies the acyclicity of $\llGCdd$.
If $S^{\la}$ stands for the set of graphs generating $\GC^\la_{marked}$, then we have
$$
\GC_{marked}^\la=\prod_{{\Ga}\in S^\la}\left(\bigotimes_{v\in V(\Ga)} C_v\right)
$$
where $C_v$ is a
\Bi
\item[(i)] one-dimensional trivial complex (generated by a black vertex) if $v$ has at least one outgoing solid edge in ${\Ga}$,
\item[(ii)] two-dimensional acyclic complex (generated by one black and one white vertex, and equipped with the differential making the white vertex black) if $v$ has no outgoing solid edges in $\Ga$.
\Ei
As any generating graph $\Ga\in S^\la$ for both possible values $\la\in \{\mathsf{t},\mathsf{or}\}$ has {\em at least one}\, vertex of type (ii), that is, at least one vertex with no outgoing solid edge, the complex $\GC_{marked}^\la$
is acyclic implying in turn acyclicity of $\llGCdd$.
\end{proof}

\newcommand{\lGCdd}{{\GC}_{d,d+1}^{\lambda}}
Let $\lGCdd$ be a subcomplex of $\llGCdd$ spanned by graphs with at least one black vertex. Then one has a short exact sequence of complexes,
$$
0\lon \lGCdd \lon \llGCdd \lon \GCd[1] \lon 0, \ \ \ \forall\ \la\in \{\mathsf{or,t}\}
$$
implying the following corollary to the Acyclicity Lemma.

\subsubsection{\bf Proposition}\label{3: Proposition on F from GCd to OGCdd} {\em For any $\la\in \{\mathsf{or,t}\}$ there is
an explicit  morphism of complexes
$$
\Ba{rccc}
F: \GCd & \lon & \lGCdd \\
   \ga & \lon & F(\ga):= (-1)^{|\ga|}\delta_{\circ\hspace{-0.5mm}-\hspace{-1mm}\bu}\ga
\Ea
$$
which is a quasi-isomorphism.
}
\begin{proof} If $\ga$ is a graph in $\lGCdd$ with all vertices white  and all edges dotted,
then $\delta \ga\equiv  \delta_{\circ\circ}\ga + \delta_{\circ\hspace{-0.5mm}-\hspace{-1mm}\bu}\ga$.
As $\delta_{\circ\circ}^2=0$ in $\OOGCdd$,  we have an identity
$$
0=\delta^2 \ga
= \delta(\delta_{\circ\hspace{-0.5mm}-\hspace{-1mm}\bu}\ga + \delta_{\circ\circ}\ga)
=\delta(\delta_{\circ\hspace{-0.5mm}-\hspace{-1mm}\bu}\ga)  +
\delta_{\circ\hspace{-0.5mm}-\hspace{-1mm}\bu}(\delta_{\circ\circ}\ga)
= (-1)^\ga\left(\delta (F(\ga)) - F(\delta_{\circ\circ} \ga)\right)
$$
which proves the claim.
\end{proof}

\sip

The complex $\lGCdd$  contain  $\lGC$ as the
subcomplex spanned by graphs with
all vertices black.

\subsubsection{\bf Lemma}\label{6: Lemma on quasi-iso OGCd to OGCdd}
 {\em For any $\la\in \{\mathsf{or,t}\}$ the inclusion $i: \lGC \hook \lGCdd$ is
 a quasi-isomorphism.}

\begin{proof} Consider a bounded above increasing filtration of $\lGC$ by the number of vertices and the number of dotted edges (so that the induced differential in the associated graded is trivial), and the following bounded above increasing filtration of $\lGCdd$
 \Beq\label{3: F_pOGCdd}
 F_{-p}\lGCdd:=\text{span of graphs $\Ga$ with}\ 2\# V_\bu(\Ga) + \#V_{\circ}(\Ga) +
 \# E'_{dot}(\Ga)\geq p
 \Eeq
 where $\# E'_{dot}(\Ga)$ is the number of dotted edges between {\em black} vertices only.
The injection $i$ respects both filtrations, hence induces an injection of the associated
graded complexes
$$
gr(i): (gr\lGC, 0) \hook (gr\lGCdd,\delta')
$$
where the differential on the right hand side complex acts only on solid edges
connecting vertices of different types (if any) by making them dotted,
$$
\delta':  \Ba{c}\resizebox{11mm}{!}{
\xy
%
 (0,0)*{\bullet}="a",
(8,0)*{\circ}="b",
\ar @{->} "a";"b" <0pt>
\endxy}\Ea
\lon
 \Ba{c}\resizebox{11mm}{!}{
\xy
%
 (0,0)*{\bullet}="a",
(8,0)*{\circ}="b",
\ar @{.} "a";"b" <0pt>
\endxy}\Ea
$$

The  Lemma is proven once we show that the map $gr(i)$ is a quasi-isomorphism which in turn is proven once we show that the direct summand $C_\circ$ of $gr\lGCdd$ spanned by graphs containing at least one white vertex is acyclic. Let us show the acyclicity of $C_\circ$.

\sip

We can assume without loss of generality that the vertices of the generating $C_\circ$ graphs are distinguished, say labelled by integers. Every such a graph contains at least one black vertex and  least one white vertex which are connected by at least one edge, and we can assume without loss of generality that their labels are 1 and 2 respectively. Then we can consider a filtration of $C_\circ$ by the number of dotted edges which connect vertices with labels not equal to 1 and 2. The associated graded complex is the tensor product of the trivial complex and the complex $C_{12}$ which controls the types of all possible edges between vertices 1 and 2. One therefore has
$$
C_{12}=\bigoplus_{k\geq 1} \odot^k C
$$
where $C$ is a 2-dimensional complex generated by the two edges,
$$
C:=\text{span}\left\langle\Ba{c}\resizebox{11mm}{!}{
\xy
(0,2)*{^{1}},
(8,2)*{^{2}},
 (0,0)*{\bullet}="a",
(8,0)*{\circ}="b",
\ar @{->} "a";"b" <0pt>
\endxy}\Ea
\  ,\
 \Ba{c}\resizebox{11mm}{!}{
\xy
(0,2)*{^{1}},
(8,2)*{^{2}},
 (0,0)*{\bullet}="a",
(8,0)*{\circ}="b",
\ar @{.} "a";"b" <0pt>
\endxy}\Ea  \right\rangle
$$
with the differential sending the solid edge into the dotted one. This complex is acyclic implying acyclicity of $C_\circ$. The Lemma is proven.
\end{proof}

\mip

The two Lemmas and the Proposition above combine together into the following theorem.

\subsubsection{\bf Theorem} {\em For any $\la\in \{\mathsf{or,t}\}$ there is a diagram of explicit quasi-isomorphisms of complexes}
$$
\GCd \lon \lGCdd \longleftarrow \lGC.
$$

\sip

This theorem implies  isomorphisms of cohomology groups in (\ref{1: GC-OGC}) and in (\ref{1: GC-TGC-SGC}).

\subsection{A useful subcomplex of the targeted complex}\label{6: subsec on GC^tdotted}  Let $\overline{\GC}_{g, d+1}^{\mathsf{t}}$ be a subcomplex of $\rTGCd$  spanned by graphs with loop number $g$. For any
graph $\Ga\in \overline{\GC}_{g, d+1}^{\mathsf{t}}$  one can write
\Beqrn
|\Ga| &=& (d+1)(\# V(\Ga) -1) + (1-d)\# E_{dot}(\Ga) + (-d) \# E_{sol}(\Ga)\\
&=& = -dg + \# V(\Ga)-1+ \# E_{dot}(\Ga)
\Eeqrn
Since every such a graph  $\Ga$ is at least trivalent, the estimation (\ref{2: V(Ga) leq 2g-2}) applies so that the degree estimation
$$
|\Ga|\leq (2-d)g -3 + \# E_{dot}(\Ga) \ \ \ \forall\Ga\in \overline{\GC}_{g, d+1}^{\mathsf{t}}.
$$
holds true.
Let
$\GC_{d+1}^{\mathsf{t\text{-}dotted}}$ be a subcomplex of $\rTGCd$ spanned by graphs having at least one dotted edge. It fits a short exact sequence
$$
0\lon \GC_{d+1}^{\mathsf{t\text{-}dotted}} \lon \rTGCd \lon  \GCD^{\mathsf{t|\geq 3}} \lon 0
$$
where the quotient complex $\GCD^{\mathsf{t|\geq 3}}$
can be identified with the subcomplex of $\dGCD^{\geq 3}$ spanned for
 by graphs with at least one target. The cohomological degree of any graph $\Ga \in \GCD^{\mathsf{t|\geq 3}}$ is bounded above:
 $$
 |\Ga|\leq (2-d)g -3.
 $$
On the other hand,  Theorem {\ref{2: lemma on vanishing H(GC)}} and the isomorphisms (\ref{1: GC-OGC}) \& (\ref{1: GC-TGC-SGC}) imply  (cf.\ \cite{Z3})
$$
H^k(\overline{\GC}_{g, d+1}^{\mathsf{t}}) \neq 0\ \ \text{only for $k$ in the range}\  (2 - d)g - 1 \leq   k \leq  (3 - d)g - 3
$$
Therefore the induced morphism of cohomology groups
$$
H^\bu(\rTGCd) \lon  H^\bu(\GCD^{\mathsf{t|\geq 3}})
$$
must be zero which proves the following

\subsubsection{\bf Proposition}\label{6: Propos on GC-dotted} {\em For any $d\in \Z$
there is a short exact sequence of cohomology groups}
$$
0\lon H^{\bu-1}(\GCD^{\mathsf{t|\geq 3}})\lon H^\bu( \GC_{d+1}^{\mathsf{t\text{-}dotted}}  )\lon H^\bu(\rTGCd) \lon 0
$$

There is a similar short exact sequence for sourced graph complexes and for oriented graph complexes. This proposition plays an important role in the  next section.



\mip

\newcommand{\fcTGCd}{\mathsf{fcGC}^{t}_{d+1}}
\newcommand{\GCvw}{\GC_{d+1}^{\vee + \wedge}}
\newcommand{\rGCvw}{\overline{\GC}_{d+1}^{\vee + \wedge}}


{\large
\section{\bf The sourced-targeted complex}
}

\subsection{An auxiliary complex} Define an auxiliary complex   $\sX_{d+1}^{\mathsf{s\cdot t}}$ by the following
short exact sequence,
$$
0\lon \GCD^{\wedge\cdot \vee} \lon \STGCd \lon \sX_{d+1}^{\mathsf{s\cdot t}} \lon 0.
$$
The complex contains a direct sum of subcomplexes spanned by graphs having at least one bivalent source
or at least one bivalent target which can be identified, respectively, with complexes $\GC_{d+1}^{\mathsf{s\text{-}dotted}}$ and $\GC_{d+1}^{\mathsf{t\text{-}dotted}}$ studied in
\S {\ref{6: subsec on GC^tdotted}} above by identifying each bivalent source (resp., target) with an
$s$-dotted edge (resp., $t$-dotted edge).
Hence there is also a short exact sequence
\Beq\label{7: SEQ on X into GC_st3}
0\lon \GC_{d+1}^{\mathsf{s\text{-}dotted}} \oplus  \GC_{d+1}^{\mathsf{t\text{-}dotted}}
\stackrel{i_1}{\lon} \sX_{d+1}^{\mathsf{s\cdot t}}\stackrel{p_1}{\lon} \GC_{d+1}^{\mathsf{s\cdot t|\geq 3}} \lon 0
\Eeq
where   the quotient complex
$\GC_{d+1}^{\mathsf{s\cdot t|\geq 3}}$
      is generated by graphs from   $\GC_{d+1}^{\mathsf{s\cdot t}}$ with all edges solid and with
      the induced differential being the standard ``vertex splitting" differential  $\delta_\bu$.


\subsection{Proposition}\label{7: Prop about H(X)} {\em There is a short exact sequence of cohomology groups}
\Beq\label{7: H(X) into GC + GC}
0\lon H^{\bu-1}(\GCD^\mathsf{{s + t|\geq 3}}) \lon  H^\bu( \sX_{d+1}^{\mathsf{s\cdot t}}) \lon
 H^\bu(\overline{\GC}_{d+1}^{\mathsf{s}}) \oplus  H^\bu(\overline{\GC}_{d+1}^{\mathsf{t}})
\lon 0
\Eeq

\begin{proof} In addition to (\ref{7: SEQ on X into GC_st3})
the complex  $\sX_{d+1}^{\mathsf{s\cdot t}}$ fits one more  short exact sequence
$$
0\lon \sX_{d+1}^{\mathsf{s\cdot t}} \stackrel{i_2}{\lon} \overline{\GC}_{d+1}^{\mathsf{s}} \oplus  \overline{\GC}_{d+1}^{\mathsf{t}}   \stackrel{p_2}{\lon} \GC_{d+1}^{\mathsf{s+ t|\geq 3}} \lon 0
$$
which stems from
following
 commutative diagram of short exact sequences of complexes,
\Beq\label{7: Comm Diagram}
\Ba{c}\resizebox{105mm}{!}{
\xymatrix{
& 0\ar[d] & 0 \ar[d] & 0\ar[d]  & \\
0\ar[r]& \GCD^{\vee\cdot \wedge} \ar[d]\ar[r] & \STGCd \ar[r]\ar[d] &
\sX_{d+1}^{\mathsf{s\cdot t}}\ar[d]^{i_2}\ar[r]  & 0\\
0\ar[r]& \GCD^\vee \oplus \GCD^\wedge \ar[d] \ar[r] & \SGCd\oplus \TGCd  \ar[r]\ar[d] &
\overline{\GC}^{\mathsf s}_{d+1} \oplus \rTGCd \ar[d]\ar[r]  & 0\\
0\ar[r]& \GCD^{\vee + \wedge}\ar[d]\ar[r] & \GCD^{\mathsf{s+t}} \ar[d]\ar[r] & \GCD^{\mathsf{s+t|\geq 3}}\ar[d]\ar[r]  & 0 \\
& 0 & 0 &  0 & \\
}
}\Ea
\Eeq
The morphism $i_2$ can be described explicitly as follows.
Any element $\Ga\in\sX_{d+1}^{\mathsf{s\cdot t}}$ is uniquely represented as a sum of three graphs
$$
\Ga=\Ga_0 + \Ga_{s} + \Ga_t
$$
where $\Ga_0$ is a graph with no bivalent vertices, $\Ga_s$ has at least one bivalent source (and hence no bivalent targets), and $\Ga_t$ has at least one bivalent target (and hence no bivalent sources).
Then the map $i_2$ is given explicitly by
$$
\Ba{rccc}
i_2: & \sX_{d+1}^{\mathsf{s\cdot t}} & \lon & \overline{\GC}^{\mathsf s}_{d+1} \oplus \rTGCd \\
& \Ga=\Ga_0+\Ga_s+\Ga_t & \lon & (\Ga_0+ \Ga_t, \Ga_0 + \Ga_s).
\Ea
$$
It respects obviously the differentials. Thanks to Proposition {\ref{6: Propos on GC-dotted}}  (and its $\mathsf{s}$-version) we conclude that  the composition $ i_2\circ i_1$ must be an epimorphism at the cohomology level. Hence the claim.
\end{proof}

\subsection{Proof of Theorem 1.1} The commutative diagram of short exact sequences (\ref{7: Comm Diagram}) implies a commutative diagram of long exact sequences

$$
\Ba{c}\resizebox{135mm}{!}{
\xymatrix{
\ldots \ar[r]& H^\bu(\GCD^\vee \oplus \GCD^\wedge) \ar[d]_{j_1} \ar[r] & H^\bu(\SGCd\oplus \TGCd)  \ar[r]\ar[d]_{j_2} &
H^\bu(\overline{\GC}^{\mathsf s}_{d+1} \oplus \rTGCd) \ar[d]_{j_3}\ar[r]  & \ldots\\
\ldots \ar[r]& H^\bu(\GCD^{\vee + \wedge})\ar[r] & H^\bu(\GCD^{\mathsf{s+t}}) \ar[r] & H^\bu(\GCD^{\mathsf{s+t|\geq 3}})\ar[r]  & \ldots \\
}
}\Ea
$$
The map $j_1$ is zero as, by Proposition {\ref{4: Prop on acyclicity of GCvee}}, the complexes
$\GCD^\vee$ and  $\GCD^\wedge$ are acyclic. The map $j_3$ is zero by the above Proposition
{\ref{7: Prop about H(X)}}. Hence the commutativity of the above diagram of cohomology groups implies that the map $j_2$ is zero as well.

\subsection{A remark on the sourced-or-targeted complex}  The complex $\GCd^{\mathsf{s+t}}$ fits the short exact sequence of complexes
$$
0\lon \GCd^{\mathsf{s+t}} \lon \dGCd \stackrel{\mathsf{p}}{\lon} \GCd^\circlearrowright \lon 0
$$
where the quotient complex $\GCd^\circlearrowright$ is generated by directed graphs $\Ga$ such that each vertex of $\Ga$ is at least trivalent and has at least one incoming and at least one outgoing edge. That complex can be identified
with the purely wheeled part $\HoLB^\circlearrowright_{p,q}(0,0)$ of the wheeled properad of strongly homotopy Lie bialgebras with $p+q+1=d$ (see \S 3.1). Composing the projection $\mathsf{p}$ with the morphism (\ref{1: GC to dGC}) we get a morphism of complexes,
$$
\mu_d: \GCd \lon \GCd^\circlearrowright,
$$
or, equivalently, a family of morphisms of their direct summands generated by graphs with loop number $g$,
$$
\mu_{g,d}: \GC_{g,d} \lon  \GC_{g,d}^\circlearrowright \ \ \ \forall \ g\geq 3.
$$
Therefore we conclude that the complex $\GCd^{\mathsf{s+t}}$ is acyclic if and only if the morphism
$\mu_{g,d}$ is a quasi-isomorphism for any $g\geq 3$. Let us test this map for the simplest non-trivial case
$g=3$, $d=2$. The left hand side of $\mu_{3,2}$ has  only one non-trivial graph in degree zero (which is a non-trivial cohomology class), the tetrahedron graph \cite{Ko1},
$$
 \Ba{c}\resizebox{13mm}{!}{
\xy
 (0,0)*{\bullet}="a",
(0,8)*{\bullet}="b",
(-7.5,-4.5)*{\bullet}="c",
(7.5,-4.5)*{\bullet}="d",
\ar @{-} "a";"b" <0pt>
\ar @{-} "a";"c" <0pt>
\ar @{-} "b";"c" <0pt>
\ar @{-} "d";"c" <0pt>
\ar @{-} "b";"d" <0pt>
\ar @{-} "d";"a" <0pt>
\endxy}
\Ea\in H^0(\GC_{3,2}).
$$
The r.h.s. of $\mu_{3,2}$ has in degree zero only two non-trivial graphs,
$$
\Ba{c}\resizebox{15mm}{!}{
\xy
 (0,0)*{\bullet}="a",
(0,8)*{\bullet}="b",
(-7.5,-4.5)*{\bullet}="c",
(7.5,-4.5)*{\bullet}="d",
\ar @{<-} "a";"b" <0pt>
\ar @{<-} "a";"c" <0pt>
\ar @{<-} "b";"c" <0pt>
\ar @{->} "d";"c" <0pt>
\ar @{->} "b";"d" <0pt>
\ar @{<-} "d";"a" <0pt>
\endxy}
\Ea \ \ \ \& \ \ \
 \Ba{c}\resizebox{16mm}{!}{
\xy
 (-5,8)*{\bullet}="a1",
  (5,8)*{\bullet}="a2",
(-5,-4.5)*{\bullet}="c",
(5,-4.5)*{\bu}="d",
\ar@/^-0.8pc/@{->}"a1";"c" <0pt>
\ar@/^0.8pc/@{<-}"a1";"c" <0pt>
\ar@/^-0.8pc/@{->}"a2";"d" <0pt>
\ar@/^0.8pc/@{<-}"a2";"d" <0pt>
\ar @{<-} "a1";"a2" <0pt>
\ar @{<-} "d";"c" <0pt>
\endxy}
\Ea
$$
which are cocycles, of course. However there is also a unique $g=3$ degree $-1$ graph in $\GC_{3,2}^\circlearrowright$
$$
 \Ba{c}\resizebox{14mm}{!}{
\xy
 (0,8)*{\bullet}="a",
(-7.5,-4.5)*{\bullet}="c",
(7.5,-4.5)*{\circ}="d",
\ar@/^-0.8pc/@{->}"a";"c" <0pt>
\ar@/^0.8pc/@{<-}"a";"d" <0pt>
\ar @{<-} "a";"c" <0pt>
\ar @{<-} "d";"c" <0pt>
\ar @{->} "a";"d" <0pt>
\endxy}
\Ea
$$
whose differential makes the above two cocyles equivalent at the cohomology level. Thus the induced map at the cohomology level $H^0(\GC_{3,2}) \rar H^0(\GC^\circlearrowright_{3,2})$ is an isomorphism.

\def\cprime{$'$}

\end{document}